\documentclass[amsppt,11pt]{amsart}
\usepackage{curves}
\usepackage{amsmath}

\newtheorem {theorem}{Theorem}[section]

\newtheorem {proposition}[theorem]{Proposition}
\newtheorem {lemma}[theorem]{Lemma}
\newtheorem {corollary}[theorem]{Corollary}

\newcounter{conjecture}
\newcounter{remark}

\newcommand{\eqnsection}{
   \renewcommand{\theequation}{\thesection.\arabic{equation}}
   \makeatletter
   \csname @addtoreset\endcsname{equation}{section}
   \makeatother}

\newcommand{\be}{{\begin{equation}}}
\newcommand{\ee}{{\end{equation}}}
\def \bt{\begin{theorem}}
\def \et{\end{theorem}}
\def \bea{\begin{eqnarray}}
\def \eea{\end{eqnarray}}
\def \bas{\begin{eqnarray*}}
\def \eas{\end{eqnarray*}}

\def \al{\alpha}
\def \bb{\beta}
\def \ga{\gamma}
\def \Ga{\Gamma}
\def \de{\delta}

\def \ep{\epsilon}
\def \epon{\ep_1}

\newcommand{\eps}{\varepsilon}
\newcommand{\bart}{{\bar{\theta}}}

\def \la{\lambda}
\def \La{\Lambda}

\def \om{\omega}

\def \si{\sigma}

\def \th{\theta}

\def \ff{\infty}
\def \wh{\widehat}
\def \wt{\widetilde}

\def \rar{\rightarrow}

\newcommand{\ls}[1]
   {\dimen0=\fontdimen6\the\font \lineskip=#1\dimen0
\advance\lineskip.5\fontdimen5\the\font \advance\lineskip-\dimen0
\lineskiplimit=.9\lineskip \baselineskip=\lineskip
\advance\baselineskip\dimen0 \normallineskip\lineskip
\normallineskiplimit\lineskiplimit \normalbaselineskip\baselineskip
\ignorespaces }
\newcommand{\req}[1]{(\ref{#1})}

\def \Z{{\bf Z}}

\def \AA{{\mathcal A}}
\def \BB{{\mathcal B}}
\def \CC{{\mathcal C}}
\def \DD{{\mathcal D}}

\def \FF{{\mathcal F}}

\def \HH{{\mathcal H}}
\def \II{{\mathcal I}}
\def \JJ{{\mathcal J}}

\def \UU{{\mathcal U}}
\def \VV{{\mathcal V}}

\def \({\left(}
\def \){\right)}

\def \lc{\left\{}
\def \rc{\right\}}

\def \nn{\nonumber}

\def \bc{\begin{center} }
\def \ec{\end{center} }

\begin{document}

\eqnsection
\newcommand{\Ini}{{I_{n,i}}}
\newcommand{\reals}{{I\!\!R}}
\newcommand{\F}{{\mathcal F}}
\newcommand{\D}{{\mathcal D}}
\newcommand{\Fn}{{{\mathcal F}_n}}
\newcommand{\Gn}{{{\mathcal G}_n}}
\newcommand{\Hn}{{{\mathcal H}_n}}
\newcommand{\Fp}{{{\mathcal F}^p}}
\newcommand{\Gp}{{{\mathcal G}^p}}
\newcommand{\PPP}{{\mathbf P}}
\newcommand{\EEE}{{\mathbf E}}
\newcommand{\Pop}{{P\otimes \PPP}}
\newcommand{\hm}{\HH^\varphi}
\newcommand{\nuw}{{\nu^W}}
\newcommand{\ths}{{\theta^*}}
\newcommand{\beq}[1]{\begin{equation}\label{#1}}
\newcommand{\eeq}{\end{equation}}
\newcommand{\integers}{{\rm I\!N}}
\newcommand{\E}{{\Bbb E}}
\newcommand{\te}{{\tilde{\delta}}}
\newcommand{\tI}{{\tilde{I}}}
\newcommand{\loge}{{\log(1/\ep)}}
\newcommand{\logen}{{\log(1/\ep_n)}}
\newcommand{\epn}{{\ep_n}}
\def\var{{\rm Var}}
\def\cov{{\rm Cov}}
\def\one{{\bf 1}}
\def\leb{{\mathcal L}eb}
\def\Ho{{\mbox{\sf H\"older}}}  
\def\thi{{\mbox{\sf Thick}}}
\def\cthi{{\mbox{\sf CThick}}}
\newcommand{\ffrac}[2]
  {\left( \frac{#1}{#2} \right)}
\newcommand{\calF}{{\mathcal F}}
\newcommand{\dfn}{\stackrel{\triangle}{=}}
\newcommand{\beqn}[1]{\begin{eqnarray}\label{#1}}
\newcommand{\eeqn}{\end{eqnarray}}
\newcommand{\oo}{\overline}
\newcommand{\uu}{\underline}
\newcommand{\bfcdot}{{\mbox{\boldmath$\cdot$}}}
\newcommand{\Var}{{\rm \,Var\,}}
\def\squarebox#1{\hbox to #1{\hfill\vbox to #1{\vfill}}}
\renewcommand{\qed}{\hspace*{\fill}
            \vbox{\hrule\hbox{\vrule\squarebox{.667em}\vrule}\hrule}\smallskip}
\newcommand{\supp}{\mbox{\rm supp}}
\newcommand{\half}{\frac{1}{2}\:}
\newcommand{\beaa}{\begin{eqnarray*}}
\newcommand{\eeaa}{\end{eqnarray*}}
\newcommand{\calK}{{\mathcal K}}
\def\dimm{{\overline{{\rm dim}}_{_{\rm M}}}}
\def\dimp{\dim_{_{\rm P}}}
\def\htaum{{\hat\tau}_m}
\def\htaumk{{\hat\tau}_{m,k}}
\def\htaumkj{{\hat\tau}_{m,k,j}}

\bibliographystyle{amsplain}

\title[Thick Points for Intersections of Planar Sample Paths]
{Thick Points for Intersections of Planar Sample Paths}

\author[Amir Dembo\,\, Yuval Peres\,\, Jay Rosen\,\, Ofer Zeitouni]
{Amir Dembo$^*$\,\, Yuval Peres$^\dagger$\,\, Jay Rosen$^\ddagger$\,\,
Ofer Zeitouni$^\S$}

\date{May 5, 2001.
\newline\indent
$^*$Research  partially supported by NSF grant \#DMS-0072331.
\newline\indent
$^\dagger$Research partially supported by NSF grant \#DMS-9803597.
\newline\indent
$^\ddagger$Research supported, in part, by grants from the NSF and
from PSC-CUNY.
\newline\indent
$^\S$
The research of all authors was supported, in part, by a US-Israel BSF grant.
\newline\indent
{\bf Key words} thick points, intersection local time, multi-fractal analysis,
stable process.
\newline\indent
{\bf AMS Subject classification:} Primary 60J55; Secondary 60J65,
28A80, 60G50.
}

\begin{abstract}
\noindent
Let $L_n^{X}(x)$ denote the number of visits to $x \in \Z^2$ of the simple
planar random walk $X$, up till step $n$. Let $X'$ be another simple
planar random walk independent of $X$. 
We show that for any $0<b<1/(2 \pi)$, there are
$n^{1-2\pi b+o(1)}$ points $x \in \Z^2$ for which
$L_n^{X}(x)L_n^{X'}(x)\geq b^2 (\log n)^4$. 
This is the discrete counterpart of
our main result, that for any $a<1$,
the Hausdorff dimension of the set of {\it thick intersection points}
$x$ for which $\limsup_{r \rightarrow 0}  \II(x,r)/(r^2|\log r|^4)=a^2$,
is almost surely $2-2a$. Here
$\II(x,r)$ is the projected intersection local time
measure of the disc of radius
$r$ centered at $x$ for two independent planar Brownian
motions run till time $1$.
The proofs rely on a `multi-scale refinement'
of the second moment method.
In addition, we also consider analogous problems
where we replace one of the Brownian motions 
by a transient stable process, or
replace the disc of radius $r$ 
centered at $x$ by $x+rK$ for general sets $K$.
\end{abstract}

\maketitle

\section{Introduction}

Let $X_n$ denote the simple random walk in the plane $\Z^2$ with
$$
L_n^{X}(x)=\#\{i: X_i=x, \,0 \leq i\leq n\}\,,\,\,\,
$$
the number of visits to 
$x\in \Z^2$ during the first $n$ steps of the walk and
$$
T_n^{X}= \max_{x\in \Z^2} L_n^{X}(x)\,,
$$
its maximal value.
In \cite[Theorem 1.1]{DPRZ4} we showed that
\begin{equation}
\label{eSRW-1a}
\lim_{n\to\infty} \frac{T_n^{X}}{(\log n)^2}=
\frac{1}{\pi}\,,
\,\,\,a.s.
\end{equation}
resolving a conjecture of Erd\"{o}s-Taylor.
Further, with 
$|A|$ counting the number of points in the set $A$, we
show there that  
\begin{equation}
\label{eSRW-2a}
\lim_{n\to\infty} \frac{\log
|\{x\in \Z^2:L_n^{X}(x)\geq a(\log n)^2 \}|}{\log n}=1-\pi a\,,\,\,\,a.s.
\end{equation}
for any fixed $0 < a < 1/\pi$.

It is natural to consider the situation for
intersections of two (or more) independent simple random walks (SRW)
in $\Z^2$. Thus, for two independent simple random walks $X_n,X_n'$
and $x\in \Z^2$, define
$$T_n^{X,X'}= \max_{x\in \Z^2} L_n^{X}(x) L_n^{X'}(x).$$
Here is our analogue of
(\ref{eSRW-1a})-(\ref{eSRW-2a}) for two independent random walks.
\begin{theorem}
\label{eSRW}
\begin{equation}
\label{eSRW-1}
\lim_{n\to\infty} \frac{T_n^{X,X'}}{(\log n)^4}=\frac{1}{4\pi^2}\,,
\,\,\,a.s.
\end{equation}
Further, 
for $0 < b < 1/(2\pi)$,
\begin{equation}
\label{eSRW-2}
\lim_{n\to\infty} \frac{\log
|\{x\in \Z^2:L_n^{X}(x)L_n^{X'}(x)\geq b^2(\log n)^4 \}|}{\log n}
=1-2\pi b\,,\,\,\,a.s.
\end{equation}
\end{theorem}
We now outline a heuristic leading to Theorem \ref{eSRW}.
With high probability the path $\{ X_k : k \leq n \}$ is contained in
a disc of radius $n^{1/2+o(1)}$.
The proof of  \req{eSRW-2a} suggests that,
with probability $n^{-\pi a+o(1)}$, a point $x$ in that disc
is visited by the random walk $\{X_k: k \leq n\}$
at least $a (\log n)^2$ times, i.e.
$L_n^{X}(x)\geq a (\log n)^2$. Similarly, the probability
that $L_n^{X'}(x)\geq a' (\log n)^2$ is $n^{-\pi a' +o(1)}$.
So, $L_n^{X}(x)L_n^{X'}(x)\geq b^2 (\log n)^4$
with probability $n^{-\pi (a+b^2/a)+o(1)}$ for some $a>0$. 
Taking the optimal value $a=b$ suggests the upper bound
given in \req{eSRW-2}. As the events
$\{ L_n^{X}(x) L_n^{X'}(x) \geq  b^2 (\log n)^4 \}$ for $x \in \Z^2$
are highly dependent, establishing the complementary lower bound
requires a `multi-scale refinement' of the second moment method.
Such a technique was developed in  \cite{DPRZ4}
for the planar Brownian occupation measure
$$
\mu_\bart^W(A):=\int_0^\bart 1_A(W_s) ds\,,
$$
where $A$ is a measurable subset of $\reals^2$,
$\{W_s, s \geq 0\}$ denotes a planar Brownian path and
$\bart:=\min\{s: |W_s|=1\}$. 
Using strong approximations
of Brownian motion by random walks,  \req{eSRW-2a} is the
outcome of a precise study of those points $x$
in whose neighborhood the Brownian motion spends an unusually large amount of
time, in all sufficiently small scales. We show that
\begin{equation}
\lim_{\ep\to 0}\sup_{x\in D(0,1)}
\frac{\mu_\bart^W(D(x,\ep))}{\ep^2 (\log \ep)^2}=2\label{nj1}
\end{equation}
where $D(x,r)$ is the open disc in $\reals^2$ of radius $r$
centered at $x$. Defining  then the set of {\it thick points}
$$\thi^W_a=\{x\in D(0,1): \lim_{\ep\to 0}
\frac{\mu_\bart^W(D(x,\ep))}{\ep^2 (\log \ep)^2}=a\}\,$$
for any $a\in (0,2]$, another outcome of
that study,
\cite[Theorem 1.4]{DPRZ4}, shows that
 $\dim(\thi^W_a)=2-a$,
with
$\dim$ denoting the Hausdorff dimension of a set.

Our goal in this paper is to better understand the essential
reason for such a result. We do so by extending it in
three different directions:

\noindent
\begin{itemize}
\item Considering occupation measure of sets without radial symmetry.
\item Considering the projected intersection local time
for two (or more) independent planar Brownian motions.
\item Considering also the projected intersection local
time between the recurrent Brownian path and the
sample path of a {\it transient} (stable) process with jumps.
\end{itemize}

Turning to the first of these
directions, let  $K\subseteq D(0,1)$
have area (i.e. Lebesgue measure) $|K|$, and let
$K(x,r)=x+rK$.
We next show how the results of \cite{DPRZ4} for thick points
are modified when we replace the discs $D(x,\ep)$ by $K(x,\ep)$. 
The analogous modification of our  study
concerning thick points for intersections will then
be straightforward and is left to the reader.
\begin{theorem}
\label{Kthick2}
If $| \partial K|=0$ then, 
\begin{equation}
\lim_{\eps\to 0}
\sup_{x \in \reals^2} \frac{\mu_\bart^W (K(x,\eps))}{\eps^2
\left(\log \frac{1}{\eps}\right)^2}=
\sup_{|x|<1} \limsup_{\eps\to0} \frac{\mu_\bart^W (K(x,\eps))}{\eps^2
\left(\log \frac{1}{\eps}\right)^2} = 2|K|/\pi
\hspace{.3in}a.s.\label{1kth}
\end{equation}
\end{theorem}
and for any $a\leq 2|K|/\pi$,
\begin{equation}
\dim\Big\{x\in \reals^2:\;
\lim_{\eps\to0} \frac{\mu_\bart^W (K(x,\eps))}{\eps^2
\left(\log \frac{1}{\eps}\right)^2}= a\Big\}= 2-a\pi/|K|
\hspace{.3in}a.s.
\label{01kth}
\end{equation}
Here one may replace $\bart$ by any
deterministic $0<T<\infty$.
We also note that all remarks following Theorem 1.3 in \cite{DPRZ4}
extend verbatim to the situation where the set $K$ is used 
instead of $D(0,1)$.

\medskip
Our main example concerns those points $x$
in whose neighborhood two independent  planar Brownian motions have an
unusually large `amount' of intersections, in all sufficiently small
scales. This is
quantified by the projected intersection local time.  Let
$\{W_s\,;\,0\leq s\leq
S\}$ and $\{W'_{t}\,;\,0\leq t\leq T\}$ be two independent planar
Brownian motions started at the
origin.
For any Borel set
$A\subseteq \reals^2$ we define the projected intersection local time in
$A$ by
\begin{equation}
\II_{S,T}(A)=\lim_{\ep\rar
0}\,\,\pi\int_0^S\int_0^T1_{A}(W_s)f_\ep(W_s-W'_{t})\,ds\,dt\label{i1.1}
\end{equation}
where the factor $\pi$ is a
convenient normalization, and $f_\ep$ is any approximate $\de$-function,
i.e. we take
$f$ to be a
non-negative
continuous function supported on the unit disc with $\int
f\,dx=1$ and set
$f_\ep(x)=f(x/\ep)/\ep^2$. It is known that the limit (\ref{i1.1}) exists
a.s. and in all
$L^p$ spaces, and that
$\II_{S,T}(\cdot)$ is a measure
supported on $\{x\in\reals^2| \,x=W_s=W'_{t} \,\,\mbox{\rm  for some $0\leq
s\leq
S,\,0\leq t\leq T$}\}$, see
\cite[Chapter VIII, Theorem 1]{Le Gall}.
Further, Le Gall \cite{Le Gall2} shows that, with $\bart' = \inf\{t:\,|W'_t| =
1\}$, there exists a constant $0<c<\infty$ such that, almost surely, for
typical $x$
in the support of $\II_{\bart,\bart'}$,
$$
\limsup_{\eps\to 0}
\frac{\II_{\bart,\,\bart'} (D(x,\eps))}{\eps^2
\left(\log \frac{1}{\eps} \log \log \log \frac{1}{\eps} \right)^2}=c\,.
$$
In contrast, our next result describes just how large the
projected intersection local time can be in the neighborhood of
any exceptional point.

\begin{theorem}
\label{theo-2}
\begin{equation}
\lim_{\eps\to 0}
\sup_{x \in \reals^2} \frac{\II_{\bart,\,\bart'} (D(x,\eps))}{\eps^2
\left(\log \frac{1}{\eps}\right)^4}=1\,,
\hspace{.6in}a.s.
\end{equation}
\end{theorem}

\noindent
Here also one may replace $\bart,\,\bart'$ by any
deterministic $0<S,T<\infty$.

The next theorem describes the multi-fractal structure
of the planar projected intersection local time.
\begin{theorem}
\label{theo-1}
For any $0<a\leq 1$,
\beq{01}
\dim\Big\{x\in D(0,1):\;
\lim_{\eps\to 0} \frac{\II_{\bart,\,\bart'} (D(x,\eps))}{\eps^2
\left(\log \frac{1}{\eps}\right)^4}= a^2\Big\}= 2-2a
\hspace{.3in}a.s.
\end{equation}
Also,
\beq{1}
\sup_{|x|<1} \limsup_{\eps\to0} \frac{\II_{\bart,\,\bart'} (D(x,\eps))}{\eps^2
\left(\log \frac{1}{\eps}\right)^4} =1
\hspace{.6in}a.s.
\end{equation}
\end{theorem}
Since $\II_{\bart,\,\bart'} (D(x,\ep))=0$
for any $x \notin \{ W_t \,\Big|\, 0 \leq t \leq \bart \}$ 
and $\ep$ small enough,
by the {\em uniform dimension doubling} property of
Brownian motion (see \cite{Kaufman} or
\cite[Eqn. (0.1)]{Perkins-Taylor}), 
\req{01} is 
equivalent to 
\beq{01t}
\dim\Big\{0\leq t\leq \bart:\;
\lim_{\eps\to0} \frac{\II_{\bart,\,\bart'} (D(W_t,\eps))}{\eps^2
\left(\log \frac{1}{\eps}\right)^4}= a^2\Big\}= 1-a
\hspace{.3in}a.s.
\end{equation}
We also have the following analogue of the coarse multi-fractal spectrum:
\begin{proposition}
\label{theo-3}
For all $a<1$,
$$ \lim_{\eps \to 0} \frac{\log \leb (x: \II_{\bart,\bart'} (D(x,\eps))\geq
a^2
\eps^2
(\log \eps)^4)}{\log \eps} = 2a\,, \quad a.s.$$
\end{proposition}

Our last example concerns
the intersections of planar Brownian paths
with other random fractals.
Let $\{ X_t \}$ denote the symmetric
stable process of index $0<\bb<2$ in the plane. As usual, we let
\begin{equation}
u^0(x)={c_{\bb}\over |x|^{2-\bb}}\label{10.1beta}
\end{equation}
denote the $0$-potential density for $\{ X_t \}$,
where $c_{\bb}=2^{-\bb} \pi^{-1}
\Ga({2-\bb\over 2})/\Ga({\bb\over 2})$. Let $\La=\La_{\bb}$ denote the norm of
\beq{Kdef}
K_{\bb}f(x)=\int_{D(0,1)} u^0(x-y)f(y)\,dy
\end{equation}
considered as an operator from $L^2\(D(0,1),\,dx\)$ to itself.

For any Borel set
$A\subseteq \reals^2$ we define the projected intersection local time in
$A$ for $W$ and
$X$ by
\begin{equation}
\II^{W,X}_{S,T}(A)=\lim_{\ep\rar
0}\,\,\pi\La^{-1}\int_0^S
\int_0^T1_{A}(W_s)f_\ep(W_s-X_{t})\,ds\,dt\label{i1.1beta}
\end{equation}
where $f_\ep$ is any approximate $\de$-function, and the factor
$\pi\La^{-1}$ is a
convenient normalization. As was the case for $\II_{S,T}(A)$ it can be
shown that
the limit (\ref{i1.1beta}) exists a.s. and in all $L^p$ spaces, and that
$\II^{W,X}_{S,T}(\cdot)$ is
a measure supported on $\{x\in\reals^2| \,x=W_s=X_{t}
 \,\,\mbox{\rm
for some
$0\leq s\leq S,\,0\leq t\leq T$}\}$. When $T=\ff$ we set
$\II^{W,X}_{S}(A):=\II^{W,X}_{S,\ff}(A)$.

We provide next the multi-fractal analysis of thick points
for $\II^{W,X}$.
\begin{theorem}
\label{theo-beta}
Fix $T \in (0,\infty]$. Then,
\begin{equation}
\lim_{\eps\to 0}
\sup_{x \in \reals^2} \frac{\II^{W,X}_{\bart,T} (D(x,\eps))}{\eps^\beta
\left(\log \frac{1}{\eps}\right)^{3}}=\beta^2/4\,,
\hspace{.6in}a.s.
\end{equation}
and for any $0<a\leq \beta/2$,
\beq{01beta}
\dim\Big\{x\in D(0,1):\;
\limsup_{\eps\to0} \frac{\II^{W,X}_{\bart,T} (D(x,\eps))}{\eps^\beta
\left(\log \frac{1}{\eps}\right)^{3}}= a^2\Big\}= \beta-2a
\hspace{.3in}a.s.
\end{equation}
Equivalently
\beq{01tbeta}
\dim\Big\{0\leq t\leq \bart:\;
\limsup_{\eps\to0} \frac{\II^{W,X}_{\bart,T} (
D(W_t,\eps))}{\eps^\beta
\left(\log \frac{1}{\eps}\right)^{3}}= a^2\Big\}= \beta/2-a
\hspace{.3in}a.s.
\end{equation}
Also,
\beq{1beta}
\sup_{|x|<1} \limsup_{\eps\to0}
\frac{\II^{W,X}_{\bart,T} (D(x,\eps))}{\eps^\beta
\left(\log \frac{1}{\eps}\right)^{3}} =
\beta^2/4
\hspace{.6in}a.s.
\end{equation}
\end{theorem}

In the next section we estimate the moments of
$\II_{\bart,\bart'} (D(x,\eps))$. Applying these estimates
we then provide the upper bounds in Theorems \ref{theo-2} and
\ref{theo-1}. Section \ref{sec-genlb} constructs the complementary
lower bounds, as well as the lower bounds for Theorem
\ref{theo-beta}, dealing with
$\II^{W,X} (D(x,\eps))$. Certain key lemmas for deriving
the lower bounds are stated and proved in Section \ref{sec-ot}.
Proposition 
\ref{theo-3} is proved in Section \ref{sec-coarse}. The
upper bounds for Theorem \ref{theo-beta} are derived in Section
\ref{sec-betau}. Section \ref{sec-kthick}
contains the proof of Theorem \ref{Kthick2}
and Section \ref{sec-etc} contains that of Theorem \ref{eSRW}.
Complements and open problems are provided in the last section.

\section{Intersection local time estimates
and upper bounds}\label{sec-hit}
\label{sec-upperbound}

Throughout this section, fix
$0<r_1\leq r$, let $\bart_r=\inf\{s>0:
|W_s|=r\}$, $\bart'_r=\inf\{t>0:
|W'_t|=r\}$
and define
$$\bar \II=\II_{\bart_{r},\,\bart'_{r}}(D(0,r_1))\,.$$
\begin{lemma}
\label{lem-hit}
We can find $c<\ff$ such that for all $k \geq 1$, $r_1\leq r/2\leq 1,$ and
$x_0,x'_0$ with
$|x_0|=|x'_0|=r_1$
\begin{equation}
\label{secmom}
\E^{x_0,\,x'_0}(\bar \II/r_1^2)^k\leq (k!)^2 \(\log
(r/r_1)+c\)^{2k}.
\end{equation}
\end{lemma}

\noindent
{\bf Proof of Lemma~\ref{lem-hit}:}
Let $g_r(x,y)$ denote the Green's function for $D(0,r)$, i.e. the
$0$-potential density
for planar Brownian motion killed when it first hits $ D(0,r)^c$.
It follows from (\ref{i1.1}) that
\begin{equation}
\E^{x_0,\,x'_0}({\bar \II}^k)=\pi^k
\sum_{\si,\si'}\int_{D(0,r_1)^k}\prod_{j=1}^k
g_{r}(y_{\si(j-1)},y_{\si(j)})g_{r}(y_{\si'(j-1)},y_{\si'(j)})\,dy_j
\label{i2.4}
\end{equation}
where the sum runs over all pairs of permutations $\si,\si'$ of
$\{1,\ldots,k\}$ and we
use the convention that $y_{\si(0)}=x_0,\,y_{\si'(0)}=x'_0$. From
\cite[page 242]{Ito-McKean} (note that our $g$ is twice theirs),
we have
$$
g_r(x,y)={1\over \pi}\log \({r\over |x-y|}|1-{x\bar{y}\over r^2}|\)
$$
where $\bar{y}$ denotes the complex conjugate of $y$. Thus, there exists
$c_o<\infty$ such that
\begin{equation}
{1\over \pi}\log \({r\over |x-y|}\)-c_o \leq 
g_r(x,y) \leq {1\over \pi}\log \({r\over |x-y|}\)+c_o
\label{i2.3}
\end{equation}
for all $|x|,\,|y|\leq r/2$.
Thus, after scaling in $r_1$, (\ref{i2.4}) is bounded by
\begin{eqnarray}
\pi^kr_1^{2k}\sum_{\si,\si'}\int_{D(0,1)^k}\prod_{j=1}^k
&&
\({1\over \pi}\log (r/r_1)-{1\over \pi}\log |y_{\si(j)}-y_{\si(j-1)}| +c_o\)
\nn\\
\label{i2.44j}
&& \({1\over \pi}\log (r/r_1)-{1\over \pi}\log
|y_{\si'(j)}-y_{\si'(j-1)}| +c_o\)\,dy_j     \\
\leq r_1^{2k}\sum_{\si,\si'}\sum_{j,j'=0}^k
&&
\(\log
(r/r_1)+c_o\pi\)^{j+j'}\sum_{\stackrel{|\AA|=k-j}{|\BB|=k-j'}} \nn
\end{eqnarray}
$$
\int_{D(0,1)^k}\prod_{i\in\AA}{}
| \log \(|y_{\si(i)}-y_{\si(i-1)}|\) | \prod_{l\in\BB}{} |\log \(
|y_{\si'(l)}-y_{\si'(l-1)}|\)| \,\frac{dy_1}{\pi}\cdots\,\frac{dy_k}{\pi} \nn
$$
where the last sum goes over all subsets $\AA,\BB$ of $\{1,\ldots,k\}$ of
cardinality
$k-j,k-j'$ respectively. We then bound the integral in the last line of
(\ref{i2.44j}) by
integrating successively with respect to $dy_1\cdots\,dy_k$, using
H\"older's inequality  and the bound
\beq{Ldef}
L:=1\vee \max_{1 \leq m \leq 4}\left(\sup_{x \in D(0,1)}\int_{D(0,1)}
 \Big|  \log \(|x-y|\)\Big|^m \, \frac{dy}{\pi}\right)^{1/m}<\ff
\end{equation}
for each integration, noting that the variable $y_j$ never appears in more
than four factors of $\log(\cdot)$. We thus bound the integral
in the last line of (\ref{i2.44j})
by $L^{(k-j)+(k-j')}$, (we use the fact that $L\geq 1$), so that
(\ref{i2.44j}) is bounded
by
\begin{eqnarray}
&&
r_1^{2k}\sum_{\si,\si'}\sum_{j,j'=0}^k\({}\log
(r/r_1)+c_o\pi\)^{j+j'}\sum_{\stackrel{|\AA|=k-j}{|\BB|=k-j'}}L^{(k-j)+(k-j')}
\nn
\\
&&
=r_1^{2k}(k!)^2\sum_{j,j'=0}^k\({}\log
(r/r_1)+c_o\pi\)^{j+j'}{k\choose j}{k\choose j'}L^{(k-j)+(k-j')} \nn
\\
&&
=r_1^{2k} (k!)^2 (\log(r/r_1)+L+c_o\pi)^{2k}\,,
\end{eqnarray}
and (\ref{secmom}) then
follows.
\qed

\begin{lemma}
\label{lem-hitprob}
For any $\de>0$ we can find $c, y_0<\ff$ and  $\eps_0>0$ so that for all
$\eps\leq \eps_0$
 and $y\geq y_0$
\begin{equation}
\label{i2.6}
P^{x_0,\,x'_0}(\II_{\bart_2,\,\bart'_2} (D(0,\eps))\geq y^2 \ep^2|\log
\ep|^2)\leq c
\exp (-(1-\de)2y)
\end{equation}
for all $x_0,\,x'_0$ with $|x_0|=|x'_0|=\eps$.
\end{lemma}

\noindent
{\bf Proof of Lemma~\ref{lem-hitprob}:}
Fix $\de>0$.
By (\ref{secmom}) with $r=2$ for
all
$\eps$ sufficiently small
and $k$ sufficiently large
\beaa
\qquad
\E^{x_0,\,x'_0}\(\II_{\bart_2,\,\bart'_2} (D(0,\eps))\)^k&\leq & (k!)^2 \(
\eps\log (1/\eps)(1+\de)\)^{2k}
\\
&\leq & (2k)!\({1\over 2}\eps\log (1/\eps)(1+
2
\de)\)^{2k}\nn
\eeaa
so that with $F={\sqrt{\II_{\bart_2,\,\bart'_2} (D(0,\eps))\over {1\over
4}\eps^2\log^2 (1/\eps)(1+
2
\de)^2}}$ we have for all $k$ large enough
$$
\E^{x_0,\,x'_0}(F^{2k})\leq (2k)!
$$
Hence, $\E^{x_0,\,x'_0}(\exp(\theta F)) \leq c_1$ for some
finite $c_1=c_1(\theta)$ and all
$|x_0|=|x'_0|=\eps$, $\theta \in [0,1)$.
Using this in Chebyscheff's inequality
yields
$$
P^{x_0,\,x'_0}(F\geq y)\leq ce^{-(1-\de)y},
$$
out of which (\ref{i2.6}) follows.
\qed

We next provide the required upper bounds in Theorems
\ref{theo-2} and \ref{theo-1}.
Namely, with the notation
\begin{equation}
\label{defthi}
 \mbox{\sf ThickInt}_{\geq a^2}=\Big\{x\in D(0,1):\;
\limsup_{\eps\to0} \frac{\II_{\bart,\,\bart'} (D(x,\eps))}{\eps^2
\left(\log \frac{1}{\eps}\right)^4}\geq a^2 \Big\} \, ,
\end{equation}
we will show that
for any $a \in (0,1]$,
\begin{equation}
\label{3.1}
\dim( \mbox{\sf ThickInt}_{\ge a^2})\leq 2-2a\,,\quad a.s.\,,
\end{equation}
and
\begin{equation}
\label{3.1b}
\limsup_{\eps\to0}
\sup_{|x| < 1}
\frac{
\II_{\bart,\,\bart'} (D(x,\eps))}{\eps^2
\left(\log \frac{1}{\eps}\right)^4} \le  1
\,,\quad a.s.
\end{equation}
(note that (\ref{3.1b}) provides the upper bound also for \req{1}).

Set $h(\ep)=\ep^2|\log \ep|^4$ and
$$
z (x,\ep) :=
\II_{\bart,\,\bart'} (D(x,\eps))/h(\eps).
$$

Fix
$\de >0$ small enough ($\delta<1/22$ will do), and
choose a sequence $\tilde\ep_n\downarrow 0$ as $n\rar
\ff$ in such a way that
$\tilde\ep_1 < e^{-2}$ and
\begin{equation}
h(\tilde\ep_{n+1}) = (1-\de) h(\tilde\ep_n),
\label{3.2}
\end{equation}
implying that $\tilde\ep_n$ is monotone decreasing in $n$.
Since, for $\tilde\ep_{n+1}\leq \ep\leq \tilde\ep_n$ we have
\begin{equation}
\label{3.3}
\quad \quad
z(x,\tilde\ep_n)= {h(\tilde\ep_{n+1})\over h(\tilde\ep_{n}) }
{\II_{\bart,\,\bart'}(D(x,\tilde\ep_n))\over
h(\tilde\ep_{n+1})} \geq
(1-\de) z(x,\ep) \, ,
\end{equation}
it is easy to see that
for any $a>0$,
\[
\mbox{\sf ThickInt}_{\geq a^2}
\subseteq D_{a^2}:=
\{x \in D(0,1) \,\Big|\,
\limsup_{n\rar\ff}z(x,\tilde\ep_n)\geq (1-\de)a^2\}.
\]
Let $\{ x_j : j=1,\ldots,\bar K_n\}$,
denote a maximal collection of points in $D(0,1)$
such that $\inf_{\ell \neq j} |x_\ell-x_j| \geq \delta \tilde\ep_n$.
Let
$\AA_{n}$ be the set of $1\leq j\leq \bar K_n$, such that
\begin{equation}
\label{defaa}
\II_{\bart,\,\bart'}(D(x_j,(1+\de)\tilde\ep_n)) \geq (1-2\de) a^2
h(\tilde\ep_n).
\end{equation}
Set $\tau_{x,n} = \inf\{s:\,W_s\in D(x,(1+\de)\tilde\ep_n)\}$, $\tau'_{x,n} =
\inf\{t:\,W'_t\in D(x,(1+\de)\tilde\ep_n)\}$. Applying the strong Markov
property and
then Lemma~\ref{lem-hitprob} we have
\beaa &&\qquad
\PPP(\II_{\bart,\,\bart'}(D(x,(1+\de)\tilde\ep_n)) \geq (1-2\de) a^2
h(\tilde\ep_n))
\\
&=& \E\(\tau_{x,n}<\bart,\,
\tau'_{x,n}<\bart'
\,;
 \right. \nn\\ &&\hspace{.7in}\left.\,
\PPP^{W_{\tau_{x,n}},W'_{\tau'_{x,n}}}(\II_{\bart,\,\bart
'}
(D(x,(1+\de)\tilde\ep_n))
\geq (1-2\de) a^2 h(\tilde\ep_n))\)\nn\\
&\leq& c\, \tilde\ep_n^{\,\,(1-10\delta)2a}\nn
\eeaa
for some $c=c(\delta)<\infty$, all sufficiently large $n$ and
any $x \in D(0,1)$.

Thus, for all sufficiently large $n$, any $j$ and $a>0$,
\begin{equation}
\label{3.3j}
\PPP(j\in \AA_n)\leq c \,\tilde\ep_n^{\,\,(1-10\de)2a}\,,
\end{equation}
implying that
\begin{equation}
\E|\AA_n|\leq c'\, \tilde\ep_n^{\,\,(1-10\de) 2a -2}.
\label{3.3jj}
\end{equation}

Let $\VV_{n,j}=D(x_j,\de\tilde\ep_n)$. For any $x \in
D(0,1)$ there exists $j \in \{1,\ldots,\bar K_n \}$ such that
$x \in \VV_{n,j}$,
hence
$D(x,\tilde\ep_n) \subseteq D(x_j,(1+\de)\tilde\ep_n)$. Consequently,
$\cup_{n \geq m} \cup_{j\in \AA_{n}}\VV_{n,j}$ forms a
cover of
$D_{a^2}$
by sets of maximal diameter $2\de \tilde\ep_m$.
Fix $a \in (0,2]$.
Since $\VV_{n,j}$ have diameter $2\de\tilde\ep_n$,
it follows from (\ref{3.3j}) that for $\gamma=2-(1-11\de)2a  > 0$,
\[
\E \sum_{n=m}^\infty \sum_{j\in \AA_{n} }|\VV_{n,j}|^\gamma \leq
c'\, (2\de)^\gamma \sum_{n=m}^\infty \tilde\ep_{n}^{\,\,\de 2a }
 < \infty \;.
\]
Thus, $\sum_{n=m}^\infty \sum_{j\in \AA_n}
|\VV_{n,j}|^\gamma$ is
 finite
a.s. implying that $\dim (D_{a^2})\leq \gamma$  a.s.
Taking $\de \downarrow 0$
completes the proof of the upper bound (\ref{3.1}).

Turning to prove (\ref{3.1b}), set $a=(1+\de)/(1-10\de)$ noting that
by (\ref{3.3jj})
$$ \sum_{n=1}^\infty \PPP(|\AA_n|\geq 1)
\leq \sum_{n=1}^\infty \E|\AA_n| \leq c'\sum_{n=1}^\infty
\tilde\ep_n^{\,\,2\de}< \infty\,.$$
By Borel-Cantelli, it follows that
a.s. $\AA_n$ is empty a.s. for all
$n>n_0(\omega)$ and some $n_0(\om)<\infty$. By (\ref{3.3}) we then have
$$
\sup_{\ep \leq \tilde\ep_{n_0(\om)}} \, \sup_{|x|<1}
\frac{\II_{\bart,\,\bart'} (D(x,\eps))}{\eps^2
\left(\log \frac{1}{\eps}\right)^4} \le  a^2\,,\,$$
and (\ref{3.1b}) follows by taking $\delta\downarrow 0$.
\qed

In proving the lower bound in the following sections we will need a variant of
Lemma~\ref{lem-hit} which we now discuss.  Let $\rho$ be a fixed measure on
$D(0,1)$
with $\int g_1(x,y)\,d\rho(y)$ uniformly bounded and continuous on $D(0,1)$
and let
$L^\rho_t$ denote the continuous additive functional for $\{W_s\}$ with Revuz
measure
$\rho$. We can define $L^\rho_t$ as
\begin{equation}
L^\rho_t=\lim_{\ep\rar 0}\int_0^\bart \int f_\ep(W_s-y)\,d\rho(y)
\label{i2.20}
\end{equation}
where $f_\ep$ is any approximate identity. Convergence in (\ref{i2.20})
holds a.s and in
all $L^p$ spaces. Our interest in such functionals stems from the fact, see
\cite[Section 2]{BK},
that for any $r$
$$
\II_{\bart_r,\,\bart'_r}
(A)=L^{\pi\mu^{W'}_{\bart'_r}|A}_{\bart_r}
$$
where $\mu^{W'}_{\bart'_r}(B)=\int_0^{\bart'_r} 1_B(W'_t)\,dt$ is
the occupation measure with respect to $W'$, and
for any measure $\nu$ we define the restriction $\nu|A(B)=\nu(A\cap B)$.
In the sequel, we write $a=(b\pm c)d$ to mean $|a/d-b|\leq c$.
\begin{lemma}
\label{lem-condhit}
Assume that for some $\gamma_1 >0$, $\gamma_2 >0$ and $c_1<\ff$
\begin{equation}
\sup_{
|x| \leq 2}
\frac{\rho (D(x,\eps))}{\eps^{\gamma_1}
\left(\log \frac{1}{\eps}\right)^{\gamma_2} }\leq c_1 <\ff\,,
\hspace{.4in}\forall \eps>0.\label{i2.23}
\end{equation}
We can then find $c<\ff$ such that for all $k \geq 1$,
$r_1\leq r_2\leq
r/2\leq
1/2
$,
and
$x_0$ with
$|x_0|=r_2$,
both
\begin{equation}
\label{csecmom}
\E^{x_0}(L^{\pi\rho|\,D(0,r_1)}_{\bart_{r}})^k\leq k!
\( \rho(D(0,r_1)) \log (r/r_1)
+ c r_1^{\gamma_1} (\log(1/r_1))^{\gamma_2} \)^k \,,
\end{equation}
and
\begin{equation}
\label{ctohelet}
\E^{x_0}(L^{\pi\rho|\,D(0,r_1)}_{\bart_{r}})=
\rho(D(0,r_1)) \log (r/r_2) \pm c r_1^{\gamma_1} (\log(1/r_1))^{\gamma_2} \,.
\end{equation}
\end{lemma}

\noindent
{\bf Proof of Lemma~\ref{lem-condhit}:}
We recall Kac's moment formula \cite{FP}, which can be derived easily from
(\ref{i2.20}):
\begin{equation}
\E^{x_0}(L^{\pi\rho|\,D(0,r_1)}_{\bart_{r}})^k=
k!\pi^k\int_{D(0,r_1)^k}\prod_{j=1}^k
g_{r}(y_{j-1},y_{j})\,d\rho(y_j).\label{i2.34}
\end{equation}
Thus, to prove (\ref{ctohelet}) we must show that
\begin{equation}
\label{ctohelet2}
\pi\int_{D(0,r_1)}
g_{r}(x_0,y)\,d\rho(y)=
\rho(D(0,r_1)) \log (r/r_2) \pm c r_1^{\gamma_1} (\log(1/r_1))^{\gamma_2} \,,
\end{equation}
for all $x_0$ with $|x_0|=r_2$, and (\ref{csecmom}) will follow from
(\ref{i2.34}) using
\begin{equation}
\label{ctohelet3}
\sup_{x\in D(0,r_1)}\pi\int_{D(0,r_1)}
g_{r}(x,y)\,d\rho(y)\leq
 \rho(D(0,r_1)) \log (r/r_1)
+ c r_1^{\gamma_1} (\log(1/r_1))^{\gamma_2} \, .
\end{equation}
inductively for the $\,dy_k,\,dy_{k-1},\ldots,\,dy_{2}$ integrals, and
(\ref{ctohelet2}) for
the $\,dy_{1}$ integral.

By (\ref{i2.3}) we can
bound the left hand side in the last display by
$$
\sup_{x\in D(0,r_1)}\int_{D(0,r_1)}\(\log (r/r_1)+\log \({r_1\over
|y-x|}\)+\pi c_o\)\,d\rho(y).
$$
Now $ D(0,r_1)\subseteq  D(x,e\,r_1)$ for any $x\in D(0,r_1)$, and
breaking
$D(x,e\,r_1)$
up into annuli
$D(x,e\,r_1)=\bigcup_{j=-1}^\ff
\lc D(x,e^{-j}r_1)\setminus D(x,e^{-j-1}r_1)\rc$ allows us to bound
\begin{eqnarray} &&
\int_{D(0,r_1)}
|\log \({r_1\over |x-y|}\)|\,d\rho(y)
\label{i2.36ba}
\\
&&\leq \sum_{j=-1}^\ff
\int_{D(x,e^{-j}r_1)\setminus
D(x,e^{-j-1}r_1)} |\log \({r_1\over
|x-y|}\)|\,d\rho(y)\nn\\
&&\leq  \sum_{j=-1}^\ff |j+2| \rho(D(x,e^{-j}r_1))\leq c_1
\sum_{j=-1}^\ff
|j+2|e^{-\gamma_1 j}r_1^{\gamma_1} (\log 1/r_1+|j|)^{\gamma_2}\nn
\end{eqnarray}
where we have used the upper bound (\ref{i2.23}) for the last
inequality. This immediately gives (\ref{ctohelet3}).

Turning to (\ref{ctohelet2}), 
by (\ref{i2.3})
we can write the left hand side as
\beaa
 &&\qquad
\int_{D(0,r_1)}
\(\log (r/r_2)+\log \({r_2\over |x_0-y|}\) \pm \pi c_o
\)\,d\rho(y)
\\
&&=\(\log (r/r_2) \pm \pi c_o
\)\rho(D(0,r_1))+\int_{D(0,r_1)}
\log \({r_2\over |x_0-y|}\)\,d\rho(y).\nn
\eeaa
When $r_1/r_2\leq 1/2$, we have that $|\log  {r_2\over |x_0-y|}|$ is bounded
on $D(0,r_1)$ so that (\ref{ctohelet2}) follows in that case.
 When $1/2\leq r_1/r_2\leq 1$ we can use
$$
\int_{D(0,r_1)}
|\log \({r_2\over |x_0-y|}\)|\,d\rho(y)\leq \int_{D(0,r_2)}
|\log \({r_2\over |x_0-y|}\)|\,d\rho(y)
$$
and (\ref{ctohelet2})
follows by the argument we applied in 
(\ref{i2.36ba}).
\qed

\section{ Lower bounds}\label{sec-lowerbound}
\label{sec-genlb}

 Fixing $a<1$, $c>0$ and $\delta>0$,
let 
\[
\Gamma_c = \Gamma_c(\om,\om') := \{x\in D(0,1):\;
\lim_{\eps\to 0} \frac{\II_{\bart_c,\,\bart'_c} (D(x,\eps))}{\eps^2
\left(\log \frac{1}{\eps}\right)^4}= a^2\} \;,
\]
 and  ${\mathcal E}_c := \{ \om,\om' : \dim(\Gamma_c(\om,\om') ) \geq 2 -2a
-\delta \}$.

In view of the results of Section \ref{sec-upperbound}, we will obtain Theorem
 \ref{theo-1} once we show that
$\PPP\times\PPP'({\mathcal E}_1)=1$ for any $a<1$ and $\de>0$.
Moreover, then the inequality
$$
\liminf_{\ep\to 0} \sup_{|x|<1}
\frac{\II_{\bart,\,\bart'} (D(x,\eps))}{\eps^2
\left(\log \frac{1}{\eps}\right)^4}
\geq
\sup_{|x|<1} \liminf_{\ep\to 0}
\frac{\II_{\bart,\,\bart'} (D(x,\eps))}{\eps^2
\left(\log \frac{1}{\eps}\right)^4}
$$
implies that for any $\eta>0$,
$$
\liminf_{\eps\to0}
\sup_{|x|<1}
\frac{\II_{\bart,\,\bart'} (D(x,\eps))}{\eps^2
\left(\log \frac{1}{\eps}\right)^4} \ge  1-\eta
\,,\quad a.s.
$$
In view of (\ref{3.1b}), these lower bounds establish Theorem \ref{theo-2}.

The bulk of this section and the next will be dedicated to showing that
$\PPP\times\PPP'({\mathcal E}_1)>0$. Assuming
this for the moment, let us show that this implies
$\PPP\times\PPP'({\mathcal E}_1)=1$.
With $W^c_t:=c^{-1}W_{c^2t}$ we have that
$c^2\bart(\om^c)= \inf\{c^2t:\,|c^{-1}W_{c^2t}| =1\}=\bart_c(\om)$, and
similarly for
$W'$ and hence
\begin{eqnarray}
\label{new-am1}
&&
\II^{W^c,W'{}^c}_{\bart,\,\bart'} (D(x,\eps))\\ &&=\lim_{\ga\rar
0}\,\,\pi\int_0^{\bart(\om^c)}\int_0^{\bart(\om'{}^c)}1_{\{|W^c_s-x|\leq\ep\}}
f_\ga (W^c_s-W'{}^c_{t})\,ds\,dt\nn\\
&&=\lim_{\ga\rar
0}\,\,\pi\int_0^{\bart(\om^c)}\int_0^{\bart(\om'{}^c)}1_{\{|W_{c^2s}-cx|\leq 
c\ep\}}
f_\ga (c^{-1}W_{c^2s}-c^{-1}W'_{c^2t})\,ds\,dt\nn\\
&&=\lim_{\ga\rar
0}\,\,{\pi\over
c^2}\int_0^{c^2\bart(\om^c)}\int_0^{c^2\bart(\om'{}^c)}1_{\{|W_{s}-cx|\leq
c\ep\}}
f_{c\ga} (W_{s}-W'_{t})\,ds\,dt\nn\\ &&= {1\over c^2}
\II_{\bart_c,\,\bart'_c} (D(cx,c\eps)). \nn
\end{eqnarray}
Consequently, $\Gamma_c(\om,\om')=c\Gamma_1(\om^c,\om'{}^c)$, so
Brownian scaling implies that  $p=\PPP\times\PPP'({\mathcal E}_c)$ is
independent of
$c>0$. Let
$$
{\mathcal E} := \limsup_{n \to \infty} {\mathcal E}_{n^{-1}} \;,
$$
so that $\PPP\times\PPP'({\mathcal E}) \geq p$.
Since
${\mathcal E}_c \in {\mathcal F}_{\bart_c}\times{\mathcal F'}_{\bart'_c}$
and $\bart_{n^{-1}} \downarrow 0$,
the Blumenthal $0-1$ law tells us that
$\PPP\times\PPP'({\mathcal E}) \in \{0,1\}$. Thus, $p>0$ yields
$\PPP\times\PPP'({\mathcal E})=1$. 
Since $x \mapsto \II_{\bart_b,\,\bart'_b} (D(x,\eps))$ are Borel
measurable,
it follows that $\Gamma_b$ are Borel sets
(hence analytic), and exactly as in \cite[(3.1)]{DPRZ4} we
find that the
events ${\mathcal E}_c$ are essentially increasing in $c$, {\it i.e.},
\begin{equation}
\label{eq:mono}
\forall 0<b<c \quad \PPP\times\PPP'({\mathcal E}_b \setminus {\mathcal
E}_c)=0 \,.
\end{equation}
Thus, $\PPP\times\PPP'({\mathcal E}\setminus {\mathcal E}_1) \leq
\PPP\times\PPP'(\bigcup_{n} \{
{\mathcal E}_{n^{-1}} \setminus {\mathcal E}_1 \}
) = 0$, so
that also $\PPP\times\PPP'({\mathcal E}_1)=1$.

It thus remains to show that $\PPP\times\PPP'({\mathcal
E}_1)>0$. We start by constructing a
subset of $\Gamma_1$, the Hausdorff dimension of which is
easier to bound below.
To this end we recall some notation from \cite{DPRZ4}.
Fix $a<2$,
$\ep_1=1/8$
and the squares
$S=S_1=[\epon,2\epon ]^2 \subset D(0,1)$. Note that
for all $x \in S$ and $y \in S \cup \{0\}$ both $0 \notin D(x,\epon)$ and
$0 \in D(x,1/2) \subset D(y,1) \subset D(x,2)$.
Let $\ep_k=\epon (k!)^{-3}=\epon \prod_{l=2}^k l^{-3}$.
For $x \in S$, $k \geq 2$
and $\rho>\epon$ let $N_k^x(\rho)$ denote the number of
excursions of $W_\cdot$ from $\partial D(x,\ep_{k-1})$ to
$\partial D(x,\ep_{k})$ prior to hitting $\partial D(x,\rho)$.
Set $n_k=3ak^2\log k$. We will say that a point $x\in S$ is
{\bf n-perfect} if
$$
n_k-k\leq N^x_k(1/2) \leq N^x_k(2) \leq n_k+k\,,
\hspace{.3in}\forall k=2,\ldots,n.
$$

For $n \geq 2$ we partition $S$ into
$M_n=\epon^2/(2 \ep_n)^2=(1/4) \prod_{l=1}^n l^6$ non-overlapping squares
of edge length $2 \ep_n=2 \epon/(n!)^3$, which we denote by
$S(n,i)\,;\,i=1,\ldots,M_n$ with $x_{n,i}$ denoting the center of
each $S(n,i)$. Let
$Y(n,i)\,;\,i=1,\ldots,M_n$ be the sequence of random variables defined by
\[Y(n,i)=1\hspace{.1in}\mbox{if $x_{n,i}$ is n-perfect}\]
and $Y(n,i)=0$ otherwise. Set
$q_{n,i}=
\PPP (Y(n,i)=1)=\E (Y(n,i))$. Define
$$
A_n=\bigcup_{i: Y(n,i)=1}S(n,i),
$$
and
\beq{F-def}
F=F(\om)=\bigcap_m \overline{\bigcup_{n\ge m} A_n}
:=\bigcap_m F_m
\,.
\end{equation}
Note that each $x \in F$ is the limit of a sequence $\{x_n\}$ such that
$x_n$ is $n$-perfect.
Finally, we set
\beq{Ca-def}
C_a:=\Big\{x\in S
: \lim_{\ep\to 0}
\frac{
\mu_{\bart'}^{w'} (D(x,\eps))}{\eps^2 \left(\log \ep\right)^2}= a\Big\},
\end{equation}

The next lemma
(to be proved in Section \ref{sec-ot}), shows
that $F\cap C_a \subset \Gamma_1$.
\begin{lemma}\label{lem1j}
A.s. for all $x\in F\cap C_a$
$$
\lim_{\eps\to 0} \frac{\II_{\bart,\,\bart'} (D(x,\eps))}{\eps^2
\left(\log \frac{1}{\eps}\right)^4}= a^2.
$$
\end{lemma}

To complete the proof
that $\PPP\times\PPP'({\mathcal E}_1)>0$ it thus suffices to show that
\begin{equation}
 \PPP\times\PPP'(\dim (F\cap C_a)\geq 2-2a- 2 \de)>0,
\label{m3.6j}
\end{equation}
for any $a<1$ and $\de>0$.

It was proved in
\cite[Section 3]{DPRZ4} that
$$
 \PPP'(\dim (C_a)\geq 2-a-\de)>0,
$$
for any $a<2$ and $\de>0$. The following general lemma will thus imply
(\ref{m3.6j}).

\begin{lemma}\label{lem-into1}
Let $B\subseteq S$ be a closed $\FF'$-measurable set with
\begin{equation}
\PPP'(\dim (B)\geq b)>0.\label{m3.6jk}
\end{equation}
Then, for any $a<b$ and $\de>0$
\begin{equation}
 \PPP\times \PPP'(\dim (F\cap B)\geq b-a-\de)>0.\label{m3.6jj}
\end{equation}
\end{lemma}

\noindent{\bf Proof.} Fixing $a<b$ and $\de>0$ such that $h:=b-a-\de >0$,
we establish (\ref{m3.6jj}) by finding a set $\mathcal{C}\times \mathcal{C}'$
 of positive $\PPP\times \PPP'$
probability,
such that for any $(\om,\om')\in \mathcal{C}\times \mathcal{C}'$
we can find a non-zero
random measure
$\rho_{\om,\om'}$ supported on $F(\om)\cap B(\om')$
with finite $h$-energy, where the
$h$-energy of a measure $\nu$ is defined as
$$
\mathcal{G}_h(\nu)=
\int\int |x-y|^{-h}\,d\nu(x)\,d\nu(y)
$$
(see e.g. \cite[Theorem 8.7]{Mat}).

By (\ref{m3.6jk})
we can find a set $\mathcal{C}'$ with
$\PPP'(\mathcal{C}')>0$ and a 
finite constant $c$
such that for any $\om'\in \mathcal{C}'$,
we have a probability measure $\mu_{\om'}$ supported on $B$
with
$$\mathcal{G}_{b-\de/2}(\mu_{\om'})<c\,.$$
In the sequel, we restrict attention to $\om'\in \mathcal{C}'$
without mentioning it explicitly.
The measure $\rho=\rho_{\om,\om'}$ will be constructed as a weak limit of
measures $\nu_n$, where $\nu_n=\nu_{n,\om,\om'}$ for $n \geq 2$
is the random measure supported on
$A_n\subseteq F_n$
whose density with respect to $\mu_{\om'}$ is
\[f_n(x)=\sum_{i=1}^{M_n}
q_{n,i}^{-1}1_{\{ Y(n,i)=1\}}1_{\{x\in S(n,i)\}}.\]
Note that
\begin{equation}
\E\(\nu_n(S)\)=\sum_{i=1}^{M_n} q_{n,i}^{-1}
\PPP
(Y(n,i)=1)
\mu_{\om'}(S(n,i))=\mu_{\om'}([\ep_1,2\ep_1]^2)=1.\label{m3.7}
\end{equation}
We next recall Lemma 3.2 of \cite{DPRZ4}, combining it with
the comments just below it.
In the sequel, we let $C_m$ denote generic finite constants that are
independent of $n$.
\begin{lemma}\label{momlb}
Let $l(i,j)=\min \{m\,: \,D(x_{n,i},\ep_{m})\cap
D(x_{n,j},\ep_{m})=\emptyset\} \leq n$,
with $l(i,i):=n$. Then,
there exists $\de_n \to 0$ such that for all $n \geq 2$, $i$,
$$
\qquad
q_{n,i} \geq
Q_n := \inf_{x \in S} \; \PPP(x\, \mbox{ is n-perfect}) \geq \ep_n^{a+\de_n} \,
\,,
$$
whereas for all $n$,
$$
\E(Y(n,i)Y(n,j))\leq C_0Q_n^2 \ep_{l(i,j)}^{-a-\de_{l(i,j)}}
\;.
$$
Furthermore, $Q_n \geq
cq_{n,i}$ for some $c>0$ and all $n\geq 2$ and $i$.
\end{lemma}
It follows from
Lemma
\ref{momlb} that
for all $\om' \in {\mathcal C}'$,
\begin{eqnarray}
\E\((\nu_n(S))^2\)&=&
\sum_{i,j=1}^{M_n}
q_{n,i}^{-1}q_{n,j}^{-1} \E(Y(n,i)Y(n,j))
\mu_{\om'}(S(n,i))\mu_{\om'}(S(n,j))
\nn\\
&\leq& C_1 \sum_{i,j=1}^{M_n}  \ep_{l(i,j)}^{-a-\de_{l(i,j)}}
\mu_{\om'}(S(n,i))\mu_{\om'}(S(n,j))
\nn\\
&\leq& C_2 \sum_{i,j=1}^{M_n} \int_{x\in S(n,i)}
\int_{y\in S(n,j)} \ep_{l(i,j)}^{-b+\de-\de_{l(i,j)}}
\mu_{\om'}(dx)\mu_{\om'}(dy)
\nn\\
&\leq &C_3 \int\int |x-y|^{-b+\delta/2}
\mu_{\om'}(dx)\mu_{\om'}(dy)
= C_3\mathcal{G}_{b-\de/2}(\mu_{\om'})\nn\\
&\leq & C_3 c <\ff
\,
\label{m3.8}
\end{eqnarray}
is a bounded sequence (recall that $\delta_l\to 0$).
Applying the Paley-Zygmund inequality (see \cite[page 8]{Kahane}),
(\ref{m3.7})
and (\ref{m3.8}) together guarantee that for some $r>0,\,v>0$ and
all $\om'\in \mathcal{C}'$,
\begin{equation}
\PPP(r^{-1} \geq \nu_n(S)\geq r)\geq 2v>0,\hspace{.2in}\forall n.
\label{pz.1}
\end{equation}
Similarly,
 for $h = b-a-\de \in (0,2)$, and $\om'\in \mathcal{C}'$,
\beaa
&& \E\(\mathcal{G}_h(\nu_n)\) \\
&\leq&
C_4 \sum_{i,j=1}^{M_n}
\frac{\E(Y(n,i)Y(n,j))}{q_{n,i}q_{n,j}}
\int_{S(n,i)}\int_{S(n,j)}
|x-y|^{-h}\,\mu_{\om'}(dx)\,\mu_{\om'}(dy) \nn\\
&\leq& C_5 \sum_{i,j=1}^{M_n} \ep_{l(i,j)}^{-a-\de_{l(i,j)}}
\int_{S(n,i)}\int_{S(n,j)}
|x-y|^{-h}\,\mu_{\om'}(dx)\,\mu_{\om'}(dy) \nn\\
&\leq & C_6
\mathcal{G}_{b-\de/2}(\mu_{\om'})\leq C_6 c<\ff.\nn
\eeaa
is a bounded sequence. Thus we can find $d<\ff$ such that
for all $\om'\in \mathcal{C}'$,
$$
\PPP(\mathcal{G}_h(\nu_n)\leq d)\geq 1-v>0,\hspace{.2in}\forall n.
$$
Combined with (\ref{pz.1}) this shows that
for all $\om'\in \mathcal{C}'$,
\begin{equation}
\PPP(r^{-1} \geq
\nu_n(S)\geq r,\,\mathcal{G}_h(\nu_n)\leq d)\geq v>0,\hspace{.2in}\forall
n.
\label{pz.3}
\end{equation}
Let
$\mathcal{C}_{n}=\{\om\,
:\, r^{-1} \geq \nu_n(S)\geq
r,\,\mathcal{G}_h(\nu_n)\leq d\,\, \mbox{\rm for all } \,\, \om'\in
\mathcal{C}'\}$
and set $\mathcal{C}=\limsup_n\mathcal{C}_{n}$. Then,
(\ref{pz.3})
implies that
$$
\PPP(\mathcal{C})\geq v>0.
$$
Fixing $\om\in \mathcal{C}$
there exists a subsequence $n_k\rar \ff$
such that $\om\in \mathcal{C}_{n_k}$ for all $k$
and all $\om'\in \mathcal{C}'$.
Due to the lower semi-continuity of $\mathcal{G}_h(\cdot)$,
the set of non-negative measures $\nu$ on $S$ such that
$\nu(S) \in [r,r^{-1}]$ and $\mathcal{G}_h(\nu) \leq d$ is compact with
respect to weak convergence.
Thus, for $(\om,\om') \in \mathcal{C}\times \mathcal{C}'$, the sequence
$\nu_{n_k}=\nu_{n_k,\om,\om'}$  has at least one weak limit $\rho_{\om,\om'}$
which is a finite measure supported on $F(\om)\cap B(\om')$,
having positive mass and finite $h$-energy.
This completes the proof of Lemma \ref{lem-into1}.
\qed

Let $\{X_t\}$ denote the symmetric stable process of index $\beta$,
with law denoted $\PPP'$.
Fix $T<\infty$ and for $a<\beta/2$, define
$$
C^X_a=\{x\in S: \limsup_{\epsilon \to 0}
\frac{\mu_T^{X} (D(x,\ep))}{\La \ep^\beta |\log \ep|}= a\}\,.
$$
As in the case of two Brownian motions, the lower bound in Theorem
\ref{theo-beta}
follows from two results. The first result
is the next lemma, to be proved in Section \ref{sec-ot}.
\begin{lemma}\label{lem1js}
A.s. for all $x\in F\cap C^X_a$
\begin{equation}
\limsup_{\eps\to 0} \frac{\II_{\bart,T}^{W,X} (D(x,\eps))}{
\eps^\beta
\left(\log \frac{1}{\eps}\right)^3}= a^2.\label{goalbjs}
\end{equation}
\end{lemma}
The second result says that
for any $\delta>0$
\begin{equation}
\PPP\times \PPP'(\dim(F\cap C^X_a)\geq \beta-2a
-\delta )>0\,.\label{rough}
\end{equation}
Indeed,
recall from \cite[(1.5)]{DPRZ3} that $\PPP'(\dim(C^X_a)=
\beta-a)=1$
(by shift invariance and stable scaling, \cite[(1.5)]{DPRZ3}
is valid even when restricting to $x \in S$).
We can extract a closed subset of $C^X_a$ which still has
dimension
$\beta-a$, so that (\ref{rough}) follows from Lemma \ref{lem-into1}.

As in the case of two Brownian motions, we conclude from
(\ref{goalbjs}) and (\ref{rough}) that
$\PPP\times \PPP'({\mathcal E}_1) >0$, where
${\mathcal E}_c := \{ \om,\om' : \dim(\Gamma_c(\om,\om') ) \geq
\beta -2a -\delta \}$ and
\[
\Gamma_c(\om,\om') := \{x\in D(0,c):\;
\limsup_{\eps\to 0} \frac{\II_{\bart_c,c^\beta T}^{W,X}(D(x,\eps))}{\eps^\beta
\left(\log \frac{1}{\eps}\right)^3}= a^2\} \;.
\]
With $W^c_t := c^{-1} W_{c^2 t}$ and $X^c_t := c^{-1} X_{c^\beta t}$,
it follows by Brownian and stable scaling that
$$
\II^{W^c,X^c}_{\bart,T} (D(x,\eps))
=  c^{-\beta} \II_{\bart_c, c^\beta T}^{W,X} (D(cx,c\eps))
$$
(see (\ref{new-am1}) for a similar derivation).
Consequently,
$\Gamma_c(\om,\om') = c \Gamma_1(\om^c,\om'{}^c)$,
with Brownian and stable scaling together implying that
$p=\PPP\times \PPP'({\mathcal E}_c) >0$ is independent of $c>0$.
Since
${\mathcal E}_c \in {\mathcal F}_{\bart_c}\times \sigma(X_t: 0 \leq t \leq c^\beta T)$,
and $\bart_{n^{-1}} \downarrow 0$,
the Blumenthal $0-1$ law tells us that
$\PPP\times\PPP'(\limsup_{n \to \infty} {\mathcal E}_{n^{-1}}) = 1$,
resulting with $\PPP \times \PPP'({\mathcal E}_1)=1$
(see (\ref{eq:mono}) for more details). This concludes the proof of 
the lower bound for any $T<\infty$. Next note that 
the random set $\Gamma_1$ for $T=\infty$ is the same as
$\Gamma_1$ for $T<\infty$ whenever $\inf_{s \geq T} |X_s| > 1$.
Consequently, in case $T=\infty$, we see that 
$$
\PPP \times \PPP' ({\mathcal E}_1) \geq \lim_{T \to \infty} 
\PPP' (\inf_{s \geq T} |X_s| > 1) =1 
$$
(by stable scaling and the transience of $X_t$).
\qed

\noindent
{\bf Remark:} The alert reader might ask whether it is possible
to get a statement similar to the one in Lemma \ref{lem-into1} but 
rather holding with probability one. Recall that
$$\thi^W_a=\{x\in D(0,1): \lim_{\ep\to 0}
\frac{\mu_\bart^W(D(x,\ep))}{\ep^2 (\log \ep)^2}=a\}\,.$$
A variation of our proof yields the
following statement:

\noindent
{\it Let $B\subset D(0,1)$ be a (possibly random)
closed set, independent 
of $\thi_a^W$. 
Further assume there exist
random probability measures
$\{\mu_{n,\omega'}(\cdot)\}$
such that
$\PPP'(\mathcal{G}_b(\mu_{1,\om'})<
\infty)=1$,
$\mu_{n,\omega'}(\cdot)$ possesses the same law as
$\mu_{1,\omega'}(\cdot/2^{n-1})$, and for any $n\geq 1$,
\begin{equation}
\label{pesah-star}
\PPP'(\mu_{n,\omega'}(B\cap D(0,2^{-(n-1)}))=1)
=1\,.
\end{equation}
Finally, assume that the sequence
of measures $\mu_{n,\omega'}$, viewed as  measure-valued random variables,
possesses a trivial tail. Then,
\begin{equation}
\label{in3.2new}
\PPP\times \PPP' (\dim  (
\thi_a^W
\cap B)\geq b-a)=1\,.
\end{equation}
}

\section{From excursions to intersection local times}\label{sec-ot}

\subsection{Intersection of two Brownian motions}
Recall the sets $F$, $C_a$,
introduced in \req{F-def} and \req{Ca-def}
and $h(\eps)=\eps^2 (\log \eps)^4$.
Lemma \ref{lem1j} will follow from the next two lemmas.

\begin{lemma}
\label{lem1l}
For every $\delta>0$, if $x\in F\cap C_a$ then
\begin{equation}
a^2
(1-\de)^5
 \le \liminf_{\ep\to 0}
\II_{\bart,\,\bart'} (D(x,\ep))/h(\ep).
\label{goalbl}
\end{equation}
\end{lemma}

\begin{lemma}
\label{lem1u}
For every $\delta>0$, if $x\in F\cap C_a$ then
\begin{equation}
\limsup_{\ep\to 0}
 \II_{\bart,\,\bart'} (D(x,\ep))/h(\ep)\le a^2
(1+\de)^5 .
\label{goalbu}
\end{equation}
\end{lemma}

\noindent{\bf Proof of Lemma \ref{lem1l}:}
Let $\de_k=\ep_k/k^6$ and let $\mathcal{D}_k$ be a $\de_k$-net of points in
$S$. Let
\[\ep'_k=\ep_ke^{1/k^{6}},\hspace{.2in}\ep''_{k-1}=\ep_{k-1}e^{-1/k^{6}},\]
so that
\begin{equation}
\ep'_k\geq\ep_k+\de_k,\hspace{.2in}\ep''_{k-1}\leq \ep_{k-1}-\de_k.\label{21.1}
\end{equation}
We will say that a point $x'\in \mathcal{D}_k$ is
{\bf lower k-successful} if there
are at least
$n_k-k$ excursions of $W$ from $\partial D(x',\ep''_{k-1})$ to
$\partial D(x',\ep'_k)$
 prior to
$\bar{\th}$.
Let \[\ep_{k,j}=\ep_ke^{-j/k}\,,\,j=0,1,\ldots,3k\log (k+1),\]
and let $\ep'_{k,j}=\ep_{k,j}e^{-2/k^{3}}=\ep'_ke^{-j/k}e^{-2/k^{3}-1/k^{6}}.$
We say that $x'\in \mathcal{D}_k$ is {\bf lower k,$\delta$-successful} if it is
lower k-successful and in addition,
\begin{equation}
\label{goalblsa}
a(1-\de) \ep^{'2}_{k,j}|\log \ep'_{k,j}|^2\le \mu^{W'}_{\bart'}
(D(x',\ep'_{k,j})) ,
\;\forall j=0,\ldots,3k\log (k+1).
\end{equation}

We now derive Lemma \ref{lem1l} from the
following lemma.
\begin{lemma}
\label{lem1ls}
There exists a $k_0=k_0(\de,\om,\om')$
such that for all $k\geq k_0$ and
$x'\in  \mathcal{D}_k$, if $x'$ is lower k,$\de$-successful
then
$$
a^2
( 1-\de)^4 h(\ep'_{k,j}) \le \II_{\bart,\,\bart'}
(D(x',\ep'_{k,j})) ,
\hspace{.2in}\forall j=0,1,\ldots,3k\log (k+1).
$$
\end{lemma}

To begin our derivation of Lemma \ref{lem1l}, note that if $x\in F$ then
there exists
a sequence of $n$-perfect points $x_n\to x$. Since
an $n$-perfect point is $k$-perfect for all $k\leq n$,
it follows that one may find a point $\tilde x_k$ which
is $k$-perfect and satisfies $|x-\tilde x_k|<\ep_k/k^6$.
Let $x_k\in \mathcal{D}_k$ with $|x_k-\tilde x_k|\leq
\ep_k/k^6$. Using (\ref{21.1}) and the fact that $\tilde x_k$ is
$k$-perfect we can see that $x_k\in \mathcal{D}_k$
is lower $k$-successful. (c.f.
\cite[Section 6, figure 1]{DPRZ4}).

Since $|x- x_k|\leq 2\ep_k/k^6$ we have
$$
 \mu^{W'}_{\bart'}(D(x_k,\ep'_{k,j}))\geq
\mu^{W'}_{\bart'}(D(x,\ep'_{k,j}-2\ep_k/k^6))
$$
for all $j,k$, so that if $x\in C_a$ it is easy to see that
there exists a $k_2=k_2(x,\om,\delta)$ such that for
$k>k_2$, $x_k$ is in fact
lower k,$\delta$-successful.

Applying Lemma \ref{lem1ls}
with $x'=x_k$ and using $\ep'_{k,j}+ 2\ep_k/k^6\leq \ep_{k,j}$ so that
$$
\II_{\bart,\,\bart'}(D(x_k,\ep'_{k,j}))\leq
\II_{\bart,\,\bart'}(D(x,\ep'_{k,j}+2\ep_k/k^6))\leq
\II_{\bart,\,\bart'}(D(x,\ep_{k,j}))
$$
then shows that for all $k$ sufficiently large
$$
a^2
(1- \de)^4 h(\ep_{k,j}) \le \II_{\bart,\,\bart'}
(D(x,\ep_{k,j})) ,
\hspace{.2in}\forall j=0,1,\ldots,3k\log (k+1).
$$
Now for any $\ep_{k+1}\leq\ep\leq \ep_k$, let $j$ be such that
$\ep_{k,j+1}\leq\ep\leq
\ep_{k,j}$. Then,
$$
\frac{\II_{\bart,\,\bart'} (D(x,\ep))}{h(\ep)}\geq
\frac{\II_{\bart,\,\bart'} (D(x,\ep_{k,j+1}))}{h(\ep_{k,j})}
\geq \frac{\II_{\bart,\,\bart'}
(D(x,\ep_{k,j+1}))}{h(\ep_{k,j+1})}\(1-4/k\)
\,,
$$
and this
completes the proof of Lemma \ref{lem1l}.

\noindent{\bf Proof of Lemma \ref{lem1ls}:}
 Suppose  that $x'\in\mathcal{D}_k$ is k,$\delta$-successful.
Then
there are at least
$n'_k=n_k-k$ excursions of $W$ between $\partial D(x',\ep'_k)$ and
$\partial D(x',\ep''_{k-1})$,
where $n'_k \to \infty$ as $k \to \infty$.
Let $\tau_{l,k,j}$ denote the projected intersection local time measure of
$D(x',\ep'_{k,j})\subset D(x',\ep'_k)$ accumulated while  $W$ executes its
$l$-th excursion between $\partial D(x',\ep'_k)$ and
$\partial D(x',\ep''_{k-1})$, and $W'$ runs from $0$ until $\bart'$.

Let
$$
A(x',k,j)
= \{\II_{\bart,\,\bart'} (D(x',\ep'_{k,j}))
\le a^2
(1-\de)^4 h(\ep'_{k,j})\}.
$$
Note that conditional on $W'$ the $\tau_{l,k,j}$ are i.i.d. and using $P^W$
to denote
probability with respect to $W$, i.e. conditional on $W'$ we have
\beaa
 P_{x',k,j}  &:=&\PPP^W \left(A(x',k,j),  x'
\mbox{\ is lower k,$\delta$-successful}\right)
\\ & \le  & \PPP_{x',s}^W \left( \sum_{l=1}^{n'_k} \tau_{l,k,j}\le
a^2(1-\de)^4 h(\ep'_{k,j}) \right) \,
\eeaa
where the subscript $x',s$ on $P^W$ indicates that
$x'$  satisfies (\ref{goalblsa}).
In \cite[Theorem 1.2]{DPRZ4} we show that
$$
\lim_{\eps\to 0}
\sup_{x \in \reals^2} \frac{\mu_{\bart'_r}^{W'} (D(x,\eps))}{\eps^2
\left(\log \frac{1}{\eps}\right)^2}=2\,,
\hspace{.6in}a.s.
$$
Consequently, the measure
$\rho(\cdot)=\mu^{W'}_{\bart'}(x'+\cdot)$ satisfies \req{i2.23}
for $\gamma_1=\gamma_2=2$.
We now apply Lemma
\ref{lem-condhit}
with $r_1=\ep'_{k,j}$,
$r_2=\ep'_k$, $ r=\ep''_{k-1}$ and $\rho(\cdot)$ as above.
Our condition (\ref{goalblsa}) implies that
$$
\rho(D(0,r_1)) = \mu^{W'}_{\bart'}(D(x',\ep'_{k,j}))
\geq a (1-\de) r_1^{\gamma_1} (\log(1/r_1))^{\gamma_2},
$$
and so, by (\ref{ctohelet}), for all $k \geq k_1(a,\delta)$,
\beq{new-am9}
E_{x',s}^W(\tau_{l,k,j}) =
\E^{x_0}(L^{\pi\rho|\,D(0,r_1)}_{\bart_{r}}) \geq
a (1-\de)^2 \log (r/r_2) r_1^{\gamma_1} (\log(1/r_1))^{\gamma_2} \,.
\end{equation}
Using Stirling's approximation for $\log \ep_k=\log \ep_1-3\log k!$,
it follows that
for all $k \geq k_2(a,\delta)$ and all $j=0,\ldots,3k\log(k+1)$,
\beq{new-am10}
n_k' \log (r/r_2) =
n_k' (3\log k- 2 k^{-6}) \geq a (1-\de) |\log (\ep'_{k,j})|^2
\end{equation}
Combining \req{new-am9} and \req{new-am10}, we see that
$$
a^2(1-\de)^3 h(\ep'_{k,j}) \leq n'_k E_{x',s}^W(\tau_{l,k,j}) \,.
$$
This immediately results with
$$
P_{x',k,j} \le \PPP_{x',s}^W \left(\frac{1}{n'_k} \sum_{l=1}^{n'_k}
\wh{\tau}_{l,k,j} \le 1-\de\right)\,,
$$
where $\wh{\tau}_{l,k,j} := \tau_{l,k,j}/E_{x',s}^W(\tau_{l,k,j})$.
Applying Lemma
\ref{lem-condhit}
as above,
we see that for all $k$ large
enough,
$$
\E_{x',s}^W (\wh{\tau}_{l,k,j} ) =1\,,
\qquad
\E_{x',s}^W(\wh{\tau}_{l,k,j} ^2) \leq
10 \,,
$$
so that, with $\wt{\tau}_{l,k,j} :=\wh{\tau}_{l,k,j}-
\E_{x',s}^W (\wh{\tau}_{l,k,j})$ we have
$$
P_{x',k,j} \le
\PPP_{x',s}^W \left(\frac{1}{n'_k} \sum_{l=1}^{n'_k} \wt{\tau}_{l,k,j}
\le -\de \right)\,.
$$
Since $\wt{\tau}_{l,k,j}\geq -1$, it follows that for all $0<\theta<1$,
$$\E_{x',s}^W(e^{-\theta \wt{\tau}_{l,k,j}})
\leq 1+2 \theta^2\E_{x',s}^W(\wt{\tau}_{l,k,j}^2)\leq 1+20\theta^2\leq
e^{20\theta^2}\,.$$
Taking $\theta=\de /40$,
a standard application of Chebyscheff's inequality then
shows that for some $\la=\la(a,\de)>0$,
$C_1<\infty$ and all $x' \in \DD_k$, $k$, $j$
$$
P_{x',k,j} \le  C_1 e^{-\la k^2\log k}\,.
$$
Since $|\DD_k| \le e^{C_2 k \log k}$ for some $C_2<\infty$ and all $k$,
it follows that
$$
\sum_{k=1}^\infty \sum_{j=0}^{3k\log (k+1)}
\sum_{x'\in \mathcal{D}_k} P_{x',k,j}
\le 3 C_1 \sum_{k=1}^\infty k^2
 e^{C_2 k\log k}   e^{-\la k^2\log k} <  \infty\,.
$$
The Borel-Cantelli lemma completes the proof of Lemma \ref{lem1ls}.
\qed

\noindent
{\bf Proof of Lemma \ref{lem1u}:} \
The situation here is quite similar to the lower
bound.
Let now
\[\bar{\ep}'_k=\ep_ke^{-2/k^{6}},\hspace{.2in}\bar{\ep}''_{k-1}=\ep_{k-1}e^{
1/k^{6}},\]
so that
$$
\bar{\ep}'_k\leq\ep_k-\de_k,\hspace{.2in}\bar{\ep}''_{k-1}\geq
\ep_{k-1}+\de_k.
$$
We now say that $x'\in \mathcal{D}_k$ is {\bf upper k-successful} if
there are at most $n_k+k$ excursions of $W$ from
$\partial D(x',\bar{\ep}''_{k-1})$ to $\partial D(x',\bar{\ep}'_k)$
prior to $\bar{\th}$.
We say it is {\bf upper k,$\de$-successful} if
it is upper k-successful and in addition,
\begin{equation}
\label{goalblsau}
a(1+\de)\ep^{'2}_{k,j}|\log \ep'_{k,j}|^2\ge \mu^{W'}_{\bart'}
(D(x',\ep'_{k,j})) ,
\hspace{.1in}\forall j=0,\ldots,3k\log (k+1).
\end{equation}

As above,
we can derive Lemma \ref{lem1u} from the  following lemma.
\begin{lemma}
\label{lem1lsu}
There exists a $k_0=k_0(\de,\om,\om')$
such that for all $k\geq k_0$ and
$x'\in  \mathcal{D}_k$, if $x'$ is upper k,$\de$-successful
then
$$
a^2
(1+\de)^4 h(\ep'_{k,j}) \ge \II_{\bart,\,\bart'}
(D(x',\ep'_{k,j})) ,
\hspace{.2in}\forall j=0,1,\ldots,3k\log (k+1).
$$
\end{lemma}

Since $0 \notin D(x,\epon)$ for all $x \in S$, with
$n''_k=n_k+k$, the proof of Lemma \ref{lem1lsu}, in analogy to that of
Lemma \ref{lem1ls}, comes down to bounding
$$
Q_{x',k,j} := \PPP_{x',us}^W \left(
 \sum_{l=1}^{n''_k} \tau_{l,k,j}
\ge a^2 (1+\de)^4 h(\ep'_{k,j}) \right) \;,
$$
where the subscript $x',us$ indicates that $x'$ satisfies
(\ref{goalblsau}).
We apply next Lemma \ref{lem-condhit} for
$r_1=\ep'_{k,j}$, $r_2=\bar{\ep}'_k$, $r=\bar{\ep}''_{k-1}$
and $\rho(\cdot)=\mu^{W'}_{\bart'}(x'+\cdot)$. Combining
(\ref{ctohelet}) and the condition (\ref{goalblsau}), it follows
that for all $k \geq k_3(a,\delta)$ and all $j=0,\ldots, 3k \log(k+1)$,
$$
 E_{x',us}^W(\tau_{l,k,j}) \leq
a (1+\de)^2 \log(\bar{\ep}''_{k-1}/\bar{\ep}'_k)
 (\ep'_{k,j})^2 |\log \ep'_{k,j}|^2 \;,
$$
and that
$$
n_k'' \log(\bar{\ep}''_{k-1}/\bar{\ep}'_k)
\leq a (1+\de) |\log (\ep'_{k,j})|^2 \;.
$$
Consequently, 
$$
a^2 (1+\de)^3 h(\ep'_{k,j}) \geq n_k'' E_{x',us}^W(\tau_{l,k,j}),
$$
and with 
$\log(\bar{\ep}''_{k-1}/\bar{\ep}'_k) \geq (1/3) \log (\bar{\ep}''_{k-1}/\ep'_{k,j})$,
it follows that for some $c_2>0$ (for example, $c_2 = a \de/3$ will do),
and all such $k,j$,
$$
Q_{x',k,j} \leq
\PPP_{x',us}^W \left(
\frac{1}{n''_k} \sum_{l=1}^{n''_k} \wt{\tau}_{l,k,j}
\ge c_2 \right)\leq e^{-c_2 \la n''_k}\(E_{x',us}^W\(e^{\la
\wt{\tau}_{1,k,j} } \)\)^{n''_k} \,,
$$
where now
$$
\wt{\tau}_{l,k,j} := \wh{\tau}_{l,k,j}-E_{x',us}^W(\wh{\tau}_{l,k,j})\,,\quad
\wh{\tau}_{l,k,j} := \frac{\tau_{l,k,j}}
{\log (\bar{\ep}''_{k-1}/\ep'_{k,j}) \ep^{'2}_{k,j}|\log \ep'_{k,j}|^2}.
$$
Applying Lemma \ref{lem-condhit}, it follows from
\req{i2.23} and \req{csecmom} that
for some $C<\infty$, all $\la>0$
small and $k$ large enough,
\beaa
E_{x',us}^W\(e^{\la  \wt{\tau}_{1,k,j} } \) &=&
1+\sum_{n=2}^\ff
\frac{\la^n}{n!} E_{x',us}^W\(\wt{\tau}^n_{1,k,j} \) \\
&\leq&
 1+\sum_{n=2}^\ff \frac{(2 \la)^n}{n!} E_{x',us}^W\(\wh{\tau}^n_{1,k,j}\) 
\leq 1+ C\la^2
\,.
\eeaa
The proof of Lemma \ref{lem1lsu} then follows as in the
proof of Lemma \ref{lem1ls}.\qed

\subsection{Intersection with stables}

Recall the sets $F$,
$C^X_a$ introduced in Section \ref{sec-genlb}.
 In this subsection, we take $T \in (0,\infty)$ and
$h(\ep)=\ep^\beta  |\log \eps|^3$.
 Lemma \ref{lem1js} is implied by the next two lemmas.

\begin{lemma}
\label{lem1lstab}
For every $\delta>0$, if $x\in F\cap C^X_a$ then
\begin{equation}
a^2(1-\de)^5  \le \limsup_{\ep\to 0}
\II^{W,X}_{\bart,T} (D(x,\ep))/h(\ep)
.\label{goalblstab}
\end{equation}
\end{lemma}
\begin{lemma}
\label{lem1ustab}
For every $\delta>0$, if $x\in F\cap C^X_a$ then
\begin{equation}
\limsup_{\ep\to 0}
 \II^{W,X}_{\bart,T} (D(x,\ep))/h(\ep)\le a^2 (1+\de)^5
.\label{goalbustab}
\end{equation}
\end{lemma}
\noindent
{\bf Proof of Lemma \ref{lem1lstab}:} Recall the notations
$\ep_k,\ep_{k,j}$, etc. of
Lemma \ref{lem1l}, and say now that $x'\in \mathcal{D}_k$
is {\bf lower k,$\delta$-successful}  if it is lower k-successful and
$$
\Lambda
a(1-\de) (\ep'_{k,j})^\beta |\log \ep'_{k,j}|\le \mu_T^{X}
(D(x',\ep'_{k,j})) ,
$$
for {\bf some} $j\in \{0,1,\ldots,3k\log (k+1)\}$.
In \cite[Theorem 1.1]{DPRZ3} we show that a.s.
$$
\lim_{\ep \to 0} \sup_{|x| \leq 2} \frac{\mu_T^X(D(x,\ep))}{\La \ep^\beta
|\log \ep|} = \beta < \infty  \,.
$$
Thus, following the arguments of the proof of Lemma \ref{lem1ls},
but now applying Lemma \ref{lem-condhit}
for $\rho(\cdot)=\La^{-1} \mu_T^{X}(x'+\cdot)$ with
$\gamma_1=\beta$ and $\gamma_2=1$ (instead of the scale
$\gamma_1=\gamma_2=2$ used in proving Lemma \ref{lem1ls}),
we see that there exists a $k_0=k_0(\de,\om,\om')$
such that for all $k\geq k_0$ and
$x'\in  \mathcal{D}_k$, if $x'$ is lower k,$\de$-successful
then
\begin{equation}
a^2 (1-\de)^4 h(\ep'_{k,j}) \le \II^{W,X}_{\bart,T}
(D(x',\ep'_{k,j})) ,
\label{goalblstab1}
\end{equation}
for some $j\in\{0,1,\ldots,3k\log (k+1)\}$.
If $x\in F$, then there exists a sequence of
points $x_k\in \mathcal{D}_k$ with $|x_k-x|\leq 2
\ep_k/k^6$
such that $x_k$ is lower k-successful. If further
$x\in C^X_a$, then there exist a
subsequence $k_n\to\infty$ and
$j_n\in \{0,1, \ldots,3k_n \log(k_n+1)\}$,
such that
$$
\Lambda
a(1-\de)(\ep'_{k_n,j_n})^\beta |
\log \ep'_{k_n,j_n}|\le \mu_T^{X}
(D(
x_{k_n},\ep'_{k_n,j_n}))\,.$$
Applying (\ref{goalblstab1}) and using the continuity of $h(\cdot)$,
one concludes the proof of Lemma \ref{lem1lstab}.
\qed

\noindent
{\bf Proof of Lemma \ref{lem1ustab}:}
The proof is analogous to that of Lemma \ref{lem1u},
now with $x' \in \mathcal{D}_k$, {\bf upper k,$\delta$-successful}
if it is upper k-successful, and such that
$$
\mu_T^{X} (D(x',\ep'_{k,j})) \leq
\Lambda a (1+\de) (\ep'_{k,j})^\beta |\log \ep'_{k,j}|,
\hspace{.1in}\forall j=0,\ldots,3k\log (k+1).
$$
Here again, the application of
Lemma \ref{lem-condhit} is with
$\gamma_1=\beta$ and $\gamma_2=1$, 
otherwise
using (\ref{csecmom}) and (\ref{ctohelet}) as in the proof of
Lemma \ref{lem1lsu}. The task of completing the
details is left to the reader.
\qed

\section{The coarse multi-fractal spectrum}\label{sec-coarse}

\noindent{\bf Proof of Proposition \ref{theo-3}:}
Fix $a\in (0,1)$ and let
$$ C(\eps,a^2)=\{x: \II_{\bart,\,\bart'} (D(x,\eps))\geq a^2\eps^2 (\log
\eps)^4\}\,.$$
With $\tilde{\eps}_n$ as in (\ref{3.2}), the bound (\ref{3.3jj})
yields for some $c_i=c_i(\delta)<\infty$ and any $\eta>0$,
\beaa
\PPP(\leb(C(\tilde{\eps}_n,a^2)) \geq \tilde{\eps}_n^\eta)
&\leq&
\tilde{\eps}_n^{-\eta} \E(\leb(C(\tilde{\eps}_n,a^2)) \\
&\leq&
c_1 \E|{\mathcal A}_n|\tilde{\eps}_n^{2-\eta} \leq
c_2\tilde{\eps}_n^{(1-10\delta)2a-\eta}\,.
\eeaa
The Borel-Cantelli lemma and (\ref{3.3}) then imply that
$$ \liminf_{\eps\to 0}
\frac{\log \leb (C(\eps,a^2/(1-\delta)))}{\log \eps}\geq 2a(1-10\delta)\,,
\quad a.s.$$
Taking $\delta\to 0$ then yields
the conclusion
$$ \liminf_{\eps\to 0}
\frac{\log \leb (C(\eps,a^2))}{\log \eps}\geq 2a\,, \quad a.s.$$
Turning to a complementary upper bound, fix $\delta>0$ such that
$a^2(1+\delta)^3<1$. Let $\ep_\de=\ep \de/(1+\de)$,
$C_\de=C(\eps/(1+\delta),a^2(1+\delta)^3)$ and $N(\ep)$ a (finite) maximal
set of $x_i \in C_\de$ such that $|x_i-x_j|>2\ep_\de$ for all $i \neq j$.
Note that $\{D(x_i,\ep_\de): x_i \in N(\ep) \}$ are disjoint and if
$x \in C_\de$ then $D(x,\ep_\de)\subset C(\eps,a^2)$. Therefore,
$$
\pi \ep_\de^2 |N(\ep)| \leq
\leb(
\bigcup_{x
\in  C_\de} D(x,\ep_\de)) \leq \leb(C(\eps,a^2))\,.
$$
With $d(\ep)=\log |N(\eps)|/\log (1/\ep)$, we thus see that
\begin{equation}
\liminf_{\ep\rar 0} d(\ep) \leq 2 -\limsup_{\eps\to 0}
\frac{\log \leb (C(\eps,a^2))}{\log \eps}\,.\label{coar.4}
\end{equation}
Let
\begin{equation}
\label{defcthi}
{\mbox{\sf CThickInt}}_{\geq a^2}=\{x\in
D(0,1):\;
\liminf_{\eps\to0} \frac{\II_{\bart,\,\bart'} (D(x,\eps))}{\eps^2
\left(\log \frac{1}{\eps}\right)^4}\geq a^2\}
\,,
\end{equation}
and
$$ {\mbox{\sf CThickInt}}_{\ga,\geq a^2}=\{x\in D(0,1):\;
\inf_{\ep \leq \ga} \frac{\II_{\bart,\,\bart'} (D(x,\eps))}{\eps^2
\left(\log \frac{1}{\eps}\right)^4}\geq a^2 \}.
$$
The sets ${\mbox{\sf CThickInt}}_{\ga,\geq a^2}$ are monotone
non-increasing in $\ga$ and
\begin{equation}{\mbox{\sf CThickInt}}_{\geq
a^2
(1+\delta)^4}\subseteq \bigcup_{n}{\mbox{\sf CThickInt}}_{\ga_n,\geq
a^2
(1+\delta)^3}\label{coar.3}
\end{equation}
for any $\ga_n\rar 0$.
Recall that  ${\mathcal S}_\ep:=\{ D(x_i,3 \ep_\de)\,:\,x_i\in  N(\ep) \}$,
forms a cover of $C_\de$, so \`{a} fortiori it is also a cover of
${\mbox{\sf CThickInt}}_{\eps/(1+\delta),\geq a^2(1+\delta)^3}$. Fixing
$\ep_n \downarrow 0$ it follows from
(\ref{coar.3}) that $\cup_{n \geq m} {\mathcal S}_{\ep_n}$ is a cover of
${\mbox{\sf CThickInt}}_{\geq
a^2
 (1+\de)^4}$ by sets of maximal diameter
$6 \ep_m$. Hence, the $\eta$-Hausdorff measure of
${\mbox{\sf CThickInt}}_{\geq a^2(1+\delta)^4}$ is finite for any $\eta$
such that
$$
\sum_{n=1}^\infty
|N(\ep_n)| \ep_n^\eta = \sum_{n=1}^\infty \ep_n^{\eta-d(\ep_n)} < \infty\,,
$$
that is, whenever
$\eta>\liminf_{\ep \to 0} d(\ep)$. Consequently, by (\ref{coar.4})
$$
\dim({\mbox{\sf CThickInt}}_{\geq a^2(1+\delta)^4})\leq \liminf_{\ep\rar 0}
d(\ep)
\leq 2 -\limsup_{\eps\to 0}
\frac{\log \leb (C(\eps,a^2))}{\log \eps}\,.
$$
Since the set considered in \req{01} is contained in
${\mbox{\sf CThickInt}}_{\geq a^2}$, taking $\delta\to 0$ and using \req{01}
yields that
$$\limsup_{\eps\to 0}
\frac{\log \leb (C(\eps,a^2))}{\log \eps}\leq 2a\,, \quad a.s. \,,$$
as needed to complete the proof.
\qed

\section{Proof of Theorem \ref{theo-beta}: the upper bound}\label{sec-betau}

Throughout this section, fix
$0<r_1\leq r$, let $\bart_r=\inf\{s>0:
|W_s|=r\}$,
and define
$$\bar \II^{W,X}=\II^{W,X}_{\bart_{r}}(D(0,r_1))\,.$$
\begin{lemma}
\label{lem-hitbeta}
For each $\delta>0$
we can find $c<\ff$ such that for all $k \geq 1$, $r_1\leq r/2\leq 1,$ and
$x_0,x'_0$ with
$|x_0|=|x'_0|=r_1$
\begin{equation}
\label{secmombeta}
\E^{x_0,\,x'_0}(\bar \II^{W,X}/r_1^\beta)^k\leq c (k!)^2 \((1+\delta)\log
(r/r_1)+c\)^{k}.
\end{equation}
\end{lemma}

\noindent
{\bf Proof of Lemma~\ref{lem-hitbeta}:}
It follows from (\ref{i1.1beta}) that
\bea
\label{i2.4betaa}
&& \E^{x_0,\,x'_0}(({\bar \II^{W,X}})^k) \\
&=&
k!(\pi\La^{-1})^k
\sum_{\si}\int_{D(0,r_1)^k}\prod_{j=1}^k
g_{r}(y_{\si(j-1)},y_{\si(j)})u^{0}(y_{j-1}-y_{j})\,dy_j
\nn
\eea
where the sum runs over all permutations $\si$ of
$\{1,\ldots,k\}$ and we
use the convention that $y_{\si(0)}=x_0$ and $y_0=x'_0$
(see \req{i2.4} for a similar formula in case of two Brownian motions).

Recalling (\ref{i2.3}), after scaling in $r_1$, (\ref{i2.4betaa}) can be
bounded 
above by
\beq{i2.44jbe}
(k!)^2 (\La^{-1})^k r_1^{k\beta} \sum_{l=0}^k\({}\log
(r/r_1)+c_o\pi\)^{l} \sum_{|\AA|=k-l} J_{\AA}
\end{equation}
where the sum goes over all subsets $\AA$ of
$\{1,\ldots,k\}$ of cardinality $k-l$, and
$$
J_{\AA} := \sup_{\si, |y_0|=|y_{\si(0)}|=1} \;
\int_{D(0,1)^k}\prod_{j\in\AA}
{}\Big |\log \(|y_{\si(j)}-y_{\si(j-1)}|\)
\Big | \prod_{j=1}^k u^{0}(y_{j-1}-y_{j})\,dy_j \,.
$$
Fix $p>1$ and $q=p/(p-1)$, and write
$$\bar L_q= \max_{m=1,2} \sup_{x\in D(0,1)}
\left(\int_{D(0,1)} \Big |\log \(|y-x|\)
\Big |^{mq} dy\right)^{1/mq}
\,.$$
Let $\La_p$ denotes the $L^2(D(0,1),dx)$ norm of the operator 
$K^{(p)}f(x)=\int_{D(0,1)} u^0(x-y)^p f(y)dy$ with kernel
$(u^{0})^p$. 
Note first that for some fixed $r$, we have that
$$C_p:= \sup_{|y_0|=1, ||f||_2\leq 1} |(K^{(p)})^r f(y)|<\infty\,.$$
Noting that each variable $y_j$  appears in at most 
a pair of logarithmic factors, we have
using H\"{o}lder's inequality that
$$ J_{\AA}\leq \bar L_q ^{k-l} \pi^{l/q} 
(C_p \La_p^{(k-r)}\pi^{1/2})^{1/p}\,.$$
Hence, by \req{i2.44jbe}, we have that
for some finite $C_p'$,
$$
\E^{x_0,\,x'_0}(({\bar \II^{W,X}})^k) \leq
C_p' (k!)^2 \left(\frac{\La_p^{1/p}}{\La}\right)^k r_1^{k\beta} 
\sum_{l=0}^k \({ \pi^{1/q}}\log
(r/r_1)+ \pi^{1+1/q}c_o\)^{l} {k\choose l} {\bar L_q}^{k-l}
$$
Taking now  $p$ large enough such that $\La_p^{1/p} \pi^{1/q}<(1+\delta) \La$
(which is possible since for $p>1$ small enough, $\La_p<\infty$, and using 
interpolation,
$\La_p\to_{p\to 1} \La$), we get
 (\ref{secmombeta}) with
$c=(c_o\pi^{1+1/q}+ \bar L_q)\vee C_p'$. 
\qed

The next lemma then follows by the
same arguments as in the
proof of Lemma~\ref{lem-hitprob},  now with
$F=\sqrt{4 \II^{W,X'}_{\bart_2} (D(0,\eps))/(
\eps^\beta |\log \eps|(1+ 4 \de))}$.

\begin{lemma}
\label{lem-hitprobbeta}
For any $\de>0$ we can find $c, y_0<\ff$ and  $\eps_0>0$ so that for all
$\eps\leq \eps_0$
 and $y\geq y_0$
\begin{equation}
\label{i2.6beta}
P^{x_0,\,x'_0}(\II^{W,X}_{\bart_2} (D(0,\eps))\geq y^2 \ep^\beta|\log
\ep|)\leq c\exp (-(1-\de)2y)
\end{equation}
for all $x_0,\,x'_0$ with $|x_0|=|x'_0|=\eps$.
\end{lemma}

The upper bound in Theorem \ref{theo-beta} is now
derived via the same line of reasoning used in \cite{DPRZ3}
(following \cite[Section 5]{DPRZ1}). The details are omitted.
\qed

\section{Proof of Theorem \ref{Kthick2}}\label{sec-kthick}

Throughout this section, fix
$0<r_1\leq r$, let $\bart_r=\inf\{s>0:
|W_s|=r\}$,
and define
$$\bar \mu={
\mu^W_{\bart_{r}}
(K(0,r_1))\over |K|/\pi}\,.$$
\begin{lemma}
\label{lem-hitkth}
We can find $c<\ff$ such that for all $k \geq 1$, $r_1\leq r_2\leq r/2\leq
1,$ and
$x_0$ with
$|x_0|=r_2$
\begin{equation}
\label{secmomkth}
\E^{x_0}(\bar \mu/r_1^2)^k\leq k! \(\log
(r/r_1)+c\)^{k}.
\end{equation}

Furthermore,
\begin{equation}
\label{ktoheletkth}
\E^{x_0}(\bar \mu/r_1^2)=\(\log
(r/r_2)\pm c\).
\end{equation}
\end{lemma}

\noindent
{\bf Proof of Lemma~\ref{lem-hitkth}:}
Let $g_r(x,y)$ denote the Green's function for $D(0,r)$, i.e. the
0-potential density
for planar Brownian motion killed when it first hits $ D(0,r)^c$.
We have
\begin{equation}
\E^{x_0}({\bar \mu}^k)=k!({\pi\over |K|})^k\int_{K(0,r_1)^k}\prod_{j=1}^k
g_{r}(y_{j-1},y_{j})\,dy_j\label{i2.4kth}
\end{equation}
where we
use the convention that $y_{\si(0)}=x_0$.

Thus, to prove (\ref{ktoheletkth}) we must show that
\begin{equation}
\label{ctohelet2kth}
{\pi\over |K|}\int_{K(0,r_1)}
g_{r}(x_0,y)\,dy=\(\log
(r/r_2)\pm c
\)r_1^2
\end{equation}
for all $x_0$ with $|x_0|=r_2$. We will also show that
\begin{equation}
\label{ctohelet3kth}
\sup_{x\in K(0,r_1)}{\pi\over |K|}\int_{K(0,r_1)}
g_{r}(x,y)\,dy\leq \(\log
(r/r_1)+
c
\)r_1^2.
\end{equation}
If we use this
 inductively for the $\,dy_k,\,dy_{k-1},\ldots,\,dy_{2}$ integrals in
(\ref{i2.4kth}), and then
use (\ref{ctohelet2kth}) for the $\,dy_{1}$ integral, we will establish
(\ref{secmomkth}).
It follows from (\ref{i2.3}) that
the left hand side in (\ref{ctohelet3kth}) is bounded by
\beaa
&&\sup_{x\in K(0,r_1)}{1\over |K|}\int_{K(0,r_1)}\(\log (r/r_1)+\log
\({r_1\over
|x-y|}\)+\pi c_o\)\,dy
\\
&&=\(\log (r/r_1)+\pi c_o\)r_1^2+r_1^2 \sup_{x\in K}{1\over |K|}\int_{K}
\log \({1\over |x-y|}\)\,dy.\nn
\eeaa
and (\ref{ctohelet3kth}) follows.

Turning to (\ref{ctohelet2kth}), as above we can write the left hand side as
\beaa &&\qquad
{1\over |K|}\int_{K(0,r_1)}
\(\log (r/r_2)+\log \({r_2\over |x_0-y|}\) \pm \pi
c_o\)\,dy
\\
&&=\(\log (r/r_2)
\pm \pi c_o
\)r_1^2+{1\over |K|}\int_{K(0,r_1)}
\log \({r_2\over |x_0-y|}\)\,dy.\nn
\eeaa
When $r_1/r_2\leq 1/2$, we have that $|\log  {r_2\over |x_0-y|}|$ is bounded
on $K(0,r_1)$ so that (\ref{ctohelet2kth}) follows in that case.
 When $1/2\leq r_1/r_2\leq 1$ we can use
$$
\int_{K(0,r_1)}
|\log \({r_2\over |x_0-y|}\)|\,dy\leq \int_{D(0,r_2)}
|\log \({r_2\over |x_0-y|}\)|\,dy
$$
 and (\ref{ctohelet2kth}) follows as above.

Continuing as in proof of Lemma~\ref{lem-hitprob},
now with $F=\bar{\mu}/\{r_1^2 |\log r_1| (1+\de)\}$,
this easily implies:

\begin{lemma}
\label{lem-hitprobkth}
For any $\de>0$ we can find $c, y_0<\ff$ and  $\eps_0>0$ so that for all
$\eps\leq \eps_0$
 and $y\geq y_0$
$$
P^{x_0}(\mu^W_{\bart_2} (K(0,\eps))\geq (|K|/ \pi)y \ep^2|\log
\ep|)\leq c\exp (-(1-\de)y)
$$
for all $x_0$ with $|x_0|=\eps$.
\end{lemma}

We first turn to the proof of the upper bounds.
The proof requires a slight adaptation of the technique of
\cite[Section 2]{DPRZ4}, because it is not true in general that
$\ep K\subset \ep' K$ if $\ep<\ep'$.

Let $r_n=(1-\de)^n$, and let $\{Q_{i,n}\}$ be a tiling
of $[-1,1]^2$ by
squares of side $r_n$, and let $Q_{i,n}^j=j+Q_{i,n}$
where $j=(j_1,j_2)$ and $j_i\in \{0,r_n/M,\ldots, (M-1)r_n/M\}$
with
$M=2/\de$. We first have the following
immediate corollary of Lemma \ref{lem-hitprobkth}.
\begin{corollary}
\label{cor-square}
There exists an $n_0=n_0(\om,\delta)$ such that
for all $n>n_0$,
$$ \#\{(i,j): \mu^W_{\bart}(Q_{i,n}^j)\geq \frac{a r_n^2 |\log r_n|^2}
{\pi}\}\leq r_n^{a(1-\de)-2}\,.$$
\end{corollary}

Our use of Corollary \ref{cor-square} is as follows.
By the assumption on the boundary of $K$, $|K|=
|\mbox{\rm cl}K|$, and hence
it is enough to prove the upper bounds for compact $K$.
Thus, let
$K\subset D(0,1)$ be compact. Fix $\delta>0$. Cover
$K$ by a finite number (say k) of squares $Q_i$, with
$\sum_{i=1}^k |Q_i|\leq |K|(1+\delta/2)$. Note that
$K(x,\ep)\subset \cup_{i=1}^k (x+\ep Q_i)$, and hence
$$\mu^W_{\bart}(K(x,\ep))\leq \sum_{i=1}^k \mu^W_{
\bart}
(x+\ep Q_i)\,.$$
Hence,
$$\{x:\limsup_{\ep\to 0}
\frac{\pi
\mu^W_{\bart}(K(x,\ep))}{|K| \ep^2 |\log \ep|^2}\geq a(1+\delta/2)\}
\subset
\bigcup_{i=1}^k
\{x:\limsup_{\ep\to 0}
\frac{\pi
\mu^W_{\bart}(Q_i(x,\ep))}{ |Q_i|\ep^2  |\log \ep|^2}\geq a\}\,.$$
Hence, it is enough to prove the upper bounds on the dimension
for an arbitrary square $Q$ of side $b$.
 Note however that then, with
$r_{n+1}\leq \ep b \leq r_n$, and any $x$, there exist
$i,j$ such that $Q(x,\ep)\subset Q_{i,n-1}^j$, whereas
$$\{x: \frac{\pi
\mu^W_{\bart}(Q(x,\ep))}{ |Q| \ep^2|\log \ep|^2}\geq a\}
\subset
\{\bigcup_{i,j} Q_{i,n-1}^j:
\frac{\pi
\mu^W_{\bart}(Q_{i,n-1}^j)}{ |Q| \ep^2|\log \ep|^2}\geq a\}
$$
Given Corollary
\ref{cor-square}, and using that $r_{n+1}/r_{n-1}=(1-\de)^2$,
one gets immediately that
$$\dim\{x:
\limsup_{\ep \to 0}
 \frac{\pi
\mu^W_{\bart}(Q_i(x,\ep))}{ |Q_i| \ep^2|\log \ep|^2}\geq a\}
\leq 2-a
(1-\de)^3
+2\de\,,\,\,\,a.s.,$$
which yields the required upper bound on the dimension.
The estimate (\ref{1kth}) is similarly proved.

Surprisingly, the lower bounds in Theorem \ref{Kthick2} require
no further computations.
First, we note the following slight adaptation of
\cite[Lemma 3.1]{DPRZ4}.
\begin{lemma}
\label{adapt}
Let $Q\subset D(0,1)$ be a fixed square, not necessarily centered.
Then there exist
 $\de(\ep)=\de(\ep,\om,Q) \to 0$ a.s.
such that for all $m$ and all $x\in D(0,
1/2)\setminus D(0,1/8)$, if $x$ is
$m$-perfect then
$$ a-\de(\ep)\leq \frac{\pi\mu_{\bart}^W(
Q(x,\ep))}{\ep^2(\log \ep)^2
|Q|}\leq a+\de(\ep)\,,\,\,\forall \ep\geq \ep_m\,.$$
\end{lemma}
\noindent
{\bf Proof of Lemma  \ref{adapt}:}
Simply follow the arguments of \cite[Section 6]{DPRZ4}
using Lemma \ref{lem-hitkth} wherever
\cite[Lemma 2.1]{DPRZ4} is used there.
\qed

Let
$$E=\{x\in D(0,1/2)\setminus D(0,1/8): \exists x_n\to x\,\,\,
\mbox{\rm  such that $x_n$ is $n$-perfect}\}\,.$$
Recall that
\begin{equation}
\label{adapt-1}
\dim(E)=2-a\,,\,\,\, a.s.
\end{equation}
It follows from Lemma \ref{adapt} that for any
fixed square $Q$, a.s.
$$ \limsup_{\ep\to 0} \frac{\pi
\mu_{\bart}^W(
Q(x,\ep))}{\ep^2(\log \ep)^2 |Q|}= a\,,
\quad \forall x\in E\,.
$$
Consequently, for any compact
$F\subset D(0,1)$, a.s.
\begin{equation}
\label{adapt0} \limsup_{\ep\to 0} \frac{\pi
\mu_{\bart}^W(
F(x,\ep))}{\ep^2(\log \ep)^2 |F|}\leq a\,,
\quad \forall x\in E\,.
\end{equation}
(Just cover $F$ by a finite number of squares $Q_i\subset D(0,1)$
with
$\sum_i |Q_i|
\leq |F|(1+\delta)$,
considering first the $\limsup$ as $\ep \to 0$, then
taking $\de\to 0$).
For $K$ satisfying the assumptions of the theorem
we have  from the above that
\begin{equation}
\label{adapt1}
 \limsup_{\ep\to 0} \frac{\pi
\mu_{\bart}^W(
K(x,\ep))}{\ep^2(\log \ep)^2 |K|}\leq
 \limsup_{\ep\to 0} \frac{\pi
\mu_{\bart}^W(
\mbox{\rm cl}K(x,\ep))}{\ep^2(\log \ep)^2 |K|}
\leq a\,,
\quad \forall x\in E\,.
\end{equation}
Recall that, for $F={\rm cl}(D(0,1) \setminus K)$ and all $x$,
\beaa
\mu_{\bart}^W(K(x,\ep))&= &
\mu_{\bart}^W(D(x,\ep))-
\mu_{\bart}^W(D(x,\ep)\setminus K(x,\ep))
\\ & \geq&
\mu_{\bart}^W(D(x,\ep))-
\mu_{\bart}^W(F(x,\ep))
\,,
\eeaa
whereas by \cite[Lemma 3.1]{DPRZ4}, a.s.
$$
 \lim_{\ep\to 0} \frac{\mu_{\bart}^W(D(x,\ep))}{\ep^2(\log \ep)^2} = a\,,
\quad \forall x\in E\,.
$$
Noting that $|F|=|D(0,1)\setminus K|$
 by our assumption on the boundary of $K$,
and using  \req{adapt0}, we conclude that a.s.
$$
\liminf_{\ep\to 0} \frac{
\mu_{\bart}^W(
K(x,\ep))}{\ep^2(\log \ep)^2 }
\geq a- \frac{a|D(0,1)\setminus K|}{\pi}= \frac{a|K|}{\pi}\,,
\quad \forall x\in E\,.
$$
When combined with \req{adapt1}, this implies that a.s.
$$
\lim_{\ep\to 0} \frac{\pi \mu_{\bart}^W(K(x,\ep))}{\ep^2(\log \ep)^2 |K|}
= a \;,
\quad \forall x\in E\,,
$$
so \req{adapt-1} yields the required lower bound
on the dimension of the sets in (\ref{01kth}).
\qed

\section{Proof of Theorem \ref{eSRW}}\label{sec-etc}

Recall the local times
$$
L_n(x):=   \sum_{k=0}^n {\bf 1}_{\{X_k=x\}}\,,\, \quad\quad
L'_n(x):=   \sum_{k=0}^n {\bf 1}_{\{X'_k=x\}}\,,\, \quad (x \in \Z^2)
$$
of the simple random walks $\{ X_k \}$ and $\{ X'_k\}$
and let $L_n^{X,X'}(x)=L_n^X(x) L_n^{X'}(x)$.
We next compute upper bounds on
$M_n^{X,X'}(b)=|\{x\in \Z^2:L_n^{X}(x)L_n^{X'}(x)\geq b \}|$
and $T_n^{X,X'}$.
To this end, using the Markov property and translation invariance of SRW in
$\Z^2$, as well as
the bound of
\cite[(5.11)]{DPRZ4}, we have for any $\delta>0$,
all $n \geq n_0(\de)$, $\alpha>0$ and $x \in \Z^2$,
$$
\PPP[L_n(x) \ge \al (\log n)^2] \leq
\PPP[L_n(0) \ge \al (\log n)^2] \leq n^{-(1-\delta) \pi \al}.
$$
Thus, fixing $0<\de<b$ and $K > (2b-\de)/\de$,
by the independence of $\{ X_k \}$ and $\{ X'_k \}$, we have
for $n \geq n_0(\de)$ and all $x \in \Z^2$,
\beqn{eofer-4}
\lefteqn{ \PPP(L_n^{X,X'}(x) \geq b^2 (\log n)^4 )}
 \\
&\leq& \sum_{i=1}^{K}
\PPP\Big[\frac{L_n(x)}{(\log n)^2} \geq (i-1)\de, \;
\frac{L'_n(x)}{(\log n)^2} \geq \frac{b^2}{i\de}\Big]
+ \PPP\Big[\frac{L_n(x)}{(\log n)^2} \geq K \de\Big]
\nn \\\hspace{.3in}
&\leq&
\sum_{i=1}^{K} n^{-(1-\delta) \pi ((i-1) \de + b^2/(i\de))} +
n^{-(1-\delta) \pi K \de}
\leq 2 K n^{-(1-\delta) \pi (2b-\de)}\nn
\eeqn
(as $s+b^2/s \geq 2b$ for all $s>0$).

Let $R_n = \max_{k \leq n} (|X_k| \vee |X'_k|)$, noting that
for some $c>0$ and all $n$ large,
\beq{eofer-3}
\PPP(R_n \geq n^{1/2+\de} ) \leq e^{-c n^\de} \;.
\end{equation}
Fixing $0<\de<b$ set
$\ga=1+4\de-(1-\de) \pi (2b -\de)$. Then, by
\req{eofer-4} and \req{eofer-3} we have for all $n$ large enough,
\beqn{eofer-6}
\lefteqn{\PPP(M_n^{X,X'}(b^2(\log n)^4) \geq n^\ga)} \nn \\
&\leq& \PPP(R_n \geq n^{1/2+\de} ) +
n^{-\ga} \E( M_n^{X,X'}(b^2(\log n)^4);\, R_n < n^{1/2+\de} )
\nn \\
&\leq& e^{-c n^\de} + n^{-\ga}
\sum_{|x| \leq n^{1/2+\de}}
 \PPP(L_n^{X,X'}(x) \geq b^2 (\log n)^4 )
\leq n^{-\de} \,.
\eeqn
For
$b<1/(2\pi)$, taking $\de \downarrow 0$
it follows by \req{eofer-6} and the Borel-Cantelli lemma that almost surely
\beq{eofer-5}
\limsup_{m \to \infty} \frac{\log M_{n_m}^{X,X'}(b^2(\log
n_m)^4)}{\log n_m} \leq
1-2 \pi b,
\end{equation}
on the subsequence $n_m=2^m$. By the monotonicity of
$n \mapsto \log L_{n}^{X,X'}(x)$ and $n \mapsto \log n$, one easily
checks that \req{eofer-5} holds also when replacing $n_m$ with $m$,
yielding the upper bound in \req{eSRW-2}. In case
$b>1/(2 \pi)$
we note that $\ga=\ga(\de)<0$ when $\de>0$ is small enough, so
\req{eofer-6} implies that for all $n$ large enough,
$$
\PPP(T_n^{X,X'} \geq b^2 (\log n)^4 ) = \PPP(M_n^{X,X'}(b^2(\log n)^4) \geq 1)
\leq n^{-\de} \;.
$$
Therefore, taking
$b \downarrow 1/(2 \pi)$,
it follows
by the Borel-Cantelli lemma that almost surely,
$$
\limsup_{m \to \infty} \frac{T_{n_m}^{X,X'}}{(\log n_m)^4} \leq
\frac{1}{4\pi^2} \;,
$$
on the subsequence $n_m=2^m$. The monotonicity of
$n \mapsto T_{n}^{X,X'}$ and $n \mapsto (\log n)^4$ allows us
to replace $n_m$ with $m$, leading to the upper bound of \req{eSRW-1}.

It suffices to prove the complementary lower bounds
for \req{eSRW-2},
because the lower bound in \req{eSRW-1} then
directly follows.
As in \cite{DPRZ4}, the proof uses the strong approximation results of
\cite{Ei}.

For any $A, A'\subseteq \Z$ let
\begin{equation}
 L_{A\times A'}^{X,X'}(z)=\#\{(i,j)\in A\times A': X_i=X'_j=z\}\label{ej.9}
\end{equation}
\begin{equation}
M\lc A\times A',\,b\rc =\#\{z\in \Z^2:
 L_{A\times A'}^{X,X'}(z)\geq b \}
\label{ej.9b}
\end{equation}
so that $M_n^{X,X'}(b)=M\lc [0,n]^2,\,b\rc$.
Let $n_{j,i}=j^{8i}e^{j^2}$ and
$\Delta n_{j,i}=n_{j,i}-n_{j,i-1}$. Fix $\ep>0$. We claim that there exists
some
$j_0=j_0(\om)<\ff$ a.s. such that for all $j\geq j_0$
\begin{equation}
\max_{1\leq i\leq j}M\lc [n_{j,i-1},n_{j,i}]^2 ,\,b^2 (\log \Delta
n_{j,i})^4\rc \geq (\Delta
n_{j,i})^{1-2\pi b-\ep}.\label{ej.11}
\end{equation}
Assuming this for the moment,
we see from (\ref{ej.11}) that for any $n_{j,j}\leq n\leq n_{j+1,j+1}$ with
$j$ sufficiently
large
\begin{equation}
M\lc [0,n]^2 ,\,(1-\ep)^2b^2 (\log n)^4\rc \geq n^{(1-\ep)^2(1-2\pi
b-\ep)}\label{ej.16}
\end{equation}
so that, replacing $b$ by $b/(1-\ep)$ we have
\begin{equation}
\liminf_{n\rar\ff}{ \log M_n^{X,X'}(b^2(\log n)^4)\over \log n}\geq
(1-\ep)^2(1- {2 \pi b\over 1-\ep}-\ep)\hspace{.3in}a.s.\label{ej.16c}
\end{equation}
and taking $\ep\rar 0$ completes the proof of the lower bound
of Theorem
\ref{eSRW} subject only to (\ref{ej.11}) which we now establish.

Note that with
$\Psi_{j,i}=\{|X_{n_{j,i-1}}| \,
\vee
\,|X'_{n_{j,i-1}}|\leq\sqrt{\Delta n_{j,i}}/\log
\Delta n_{j,i}\}$ we have
\begin{equation}
\PPP\(\max_{1\leq i\leq j}M\lc [n_{j,i-1},n_{j,i}]^2 ,\,b^2 (\log \Delta
n_{j,i})^4\rc \leq
(\Delta n_{j,i})^{1-2\pi b-\ep}\)\label{ej.16cn}
\end{equation}
\[\leq j\max_{1\leq i\leq j}\PPP (\Psi^c_{j,i})\]
\[+\PPP\(\bigcap_{1\leq i\leq j}\lc M\lc [n_{j,i-1},n_{j,i}]^2 ,\,b^2 (\log
\Delta n_{j,i})^4\rc \leq (\Delta n_{j,i})^{1-2\pi b-\ep}\,;\,\Psi_{j,i}\rc\)\]
We will show that this is summable in $j$ so that (\ref{ej.11}) will follow by
the Borel-Cantelli Lemma.  By \cite[Theorem 17.5]{Revesz}, (which is
essentially the
Central Limit Theorem), for
some $c>0$ and all
$j$ sufficiently large
\begin{equation}
\PPP (|X_{n_{j,i-1}}|
>\sqrt{\Delta n_{j,i}}/\log \Delta n_{j,i})\leq \PPP
({|X_{n_{j,i-1}}|\over
\sqrt{n_{j,i-1}} }>
j)\leq e^{-cj}
\label{ej.13}
\end{equation}
so that the first term on the right hand side of (\ref{ej.16cn}) is
summable in $j$.
On the other hand, since
\bas && M\lc [n_{j,i-1},n_{j,i}]^2 ,\,b^2 (\log \Delta n_{j,i})^4\rc\\
&&=M\lc [0,\Delta
n_{j,i}]^2
,\,b^2 (\log
\Delta n_{j,i})^4\rc \circ
(\theta_{n_{j,i-1}},\,\theta'_{n_{j,i-1}}),\eas
 the Markov property together with the
following lemma will bound the second term on the right hand side of
(\ref{ej.16cn}) by
$(1-p_0)^j$ which is also
summable in $j$.  Thus,
with $D_0(r)=D(0,r) \cap \Z^2$,
(\ref{ej.11}) and consequently the proof of
the theorem is reduced to the following lemma.

\begin{lemma}\label{lem-eopen}
For any $\ep>0$ we can find $ p_0>0$ such that for all $n $ sufficiently
large and all
$x,x'\in D_0(\sqrt{n}/\log n)$
\begin{equation}
\PPP^x \times {\PPP'}^{x'} \(M_n^{X,X'}(b^2(\log n)^4) \geq n^{1-2\pi
b-\ep} \) \geq
p_0>0.\label{ej.10}
\end{equation}
\end{lemma}
\qed

\noindent
{\bf Proof of Lemma \ref{lem-eopen}}. For any $B\subseteq \Z^2$, let
\begin{equation}
M\lc A\times A',\,B,\,b\rc =\#\{z\in B:
 L_{A\times A'}^{X,X'}(z)\geq b \}\label{ej.17}
\end{equation}
and note that if $\tau_B=\inf\{
i \geq 0 : X_i\in B\}$, $\tau'_B=\inf\{
i \geq 0 :
X'_i\in B\}$ we have
\begin{equation}
M\lc A\times A',\,B,\,b\rc=M\lc (A-\tau_B)\times
(A'-\tau'_B),\,B,\,b\rc\circ (\th_{\tau_B},\,\th'_{\tau'_B}).\label{ej.17a}
\end{equation}
We will use the abbreviation $\wt{b}=1-2\pi b$.
Let $\tau_{n}=\tau_{D^c_0(2\sqrt{n}/\log n)}$ and
$\tau'_{n}=\tau'_{D^c_0(2\sqrt{n}/\log
n)}$. By (\ref{ej.17a}) and the Markov property
\begin{equation}
\PPP \times {\PPP'} \(M\lc [0,n]^2 ,\,D^c_0(2\sqrt{n}/\log n),\,b^2 (\log
n)^4\rc
\geq n^{\wt{b}-\ep} \)\label{ej.19}
\end{equation}
\[\leq \E \times \E'\lc
\PPP^{X_{\tau_{n}}} \times {\PPP'}^{X'_{\tau_{n}}} \(M\lc [0,n]^2
,\,D^c_0(2\sqrt{n}/\log n),\,b^2 (\log n)^4\rc \geq n^{\wt{b}-\ep}
\)\rc\]
\begin{eqnarray} &&
=\sum_{y,y'}H_{n}(0,y)H_{n}(0,y')\nn\\
&&\hspace{.8in}
\PPP^{y} \times {\PPP'}^{y'} \(M\lc [0,n]^2
,\,D^c_0(2\sqrt{n}/\log n),\,b^2 (\log n)^4\rc \geq n^{\wt{b}-\ep}
\)\nn
\end{eqnarray}
where
$H_{n}(x,y)=\PPP^x(X_{\tau_{n}}=y)$ is a
harmonic measure. By Harnack's
inequality, \cite[Theorem 1.7.2]{lawler}, for some $C<\ff$ and all $n,y$
\begin{equation}
H_{n}(0,y)\leq C\inf_{x\in D_0(\sqrt{n}/\log n)}H_{n}(x,y)\label{ej.20}
\end{equation}
Using (\ref{ej.17a}) and the Markov
property again shows that uniformly in $x,x'\in D_0(\sqrt{n}/\log n)$
\begin{eqnarray}&&\hspace{.2in}
\PPP \times {\PPP'} \(M\lc [0,n]^2 ,\,D^c_0(2\sqrt{n}/\log n),\,b^2 (\log
n)^4\rc
\geq n^{\wt{b}-\ep} \)\label{ej.21}\\
&&
\leq C\PPP^x \times {\PPP'}^{x'} \(M\lc [0,n+\tau_{n}]\times
[0,n+\tau'_{n}],\,b^2 (\log n)^4\rc \geq n^{\wt{b}-\ep} \)\nn\\
&&
\leq C\PPP^x \times {\PPP'}^{x'} \(M\lc [0,(1+\ep)n]^2
,\,b^2 (\log n)^4\rc \geq n^{\wt{b}-\ep}
\) \nn\\
&&  \hspace{1in} +C \PPP^x(\tau_{n}>\ep n)+
C \PPP^{x'}(\tau'_{n}>\ep n).\nn
\end{eqnarray}
Since $\PPP^x(\tau_{n}>\ep n)=\PPP^x(\max_{j\leq \ep n}|X_j|<
2\sqrt{n}/\log n)$, we see
from the last line of
\cite[Theorem 17.5]{Revesz} that for $n$ sufficiently large the last line
of (\ref{ej.21}) is
negligible,  so that given  the next lemma, we obtain
Lemma \ref{lem-eopen} after some adjustment of $b$ and $\ep$. \qed

\begin{lemma}\label{lem-esecond}
For any $\ep>0$ we can find $ p_0>0$ such that for all $n $ sufficiently large
\begin{equation}
\PPP \times {\PPP'} \(M\lc [0,n]^2,\,D^c_0(2\sqrt{n}/\log n) ,\,b^2 (\log
n)^4\rc
\geq n^{1-2\pi b-
3 \ep}
\)
\geq p_0.\label{ej.18}
\end{equation}
\end{lemma}

\noindent
{\bf Proof of Lemma \ref{lem-esecond}} This lemma will be derived
 from results about Brownian motion by using strong approximation. However,
results
about $M\lc [0,n]^2,\,D^c_0(2\sqrt{n}/\log n) ,\,b^2 (\log
n)^4\rc$ are not easy to obtain directly by strong approximation since
$L_{[0,n]^2}^{X,X'}(z)$ does not correspond to a functional of Brownian
motion. Instead,
we will derive results about $L_{[0,n]^2}^{X,X'}(z)$ from results about
excursions of
random walks between concentric discs, and it is such results which can be
obtained
from our work on Brownian thick points. To this end we introduce notation
which is
meant to simplify the connection with Brownian motion.

Fix
$a<1$ and
$\de>0$ small.
Set $m_k=ak^2$ and
$k(n)=[(1/2-\de)\log n]$ so that $e^{-k(n)}\geq n^{\de-1/2}$.
Let $D_z(r)=D(z,r) \cap \Z^2$
denote the disc of radius $r$ in $\Z^2$ centered at
$z$ and define its boundary
$\partial D_z(r) = \{ z' \notin D_z(r) : |z'-y|=1 \mbox{ for
some } y \in D_z(r) \}$.
For any fixed $K<\infty$ let
\[r_n=(1+3 \de) e^{-k(n)} \sqrt{n/(2K)},\hspace{.3in}
R_n =(1-3 \de) e^{-k(n)+1} \sqrt{n/(2K)}\]
We say that $z \in \Z^2$ is
{\bf $n, \de$-admissible} if
at least $(1-2\de) m_{k(n)}$ excursions between $\partial D_z(r_n)$ and
$\partial D_z(R_n)$ are completed by
both $\{X_i \}$ and $\{X'_i \}$ for $i \leq n$.
Let
$$ I_{n, \de} := \{ z \in D_0(\sqrt{n/(2K)})\cap D^c_0(
2 \sqrt{n}/\log n) : z
\,\mbox{\rm is $n, \de$-admissible}\}\,.$$
Lemma \ref{lem-esecond} now follows from the next lemma if we take $a=2\pi
b/(1-\ep)$ and
$\ep$,
$\de$ sufficiently small.\qed

\begin{lemma}
\label{eSRW-occ} For any $\ep>0$ we can find $K_0<\ff$ and $\de_0>0$
 such that for all $K>K_0$ and $\de<\de_0$
\beq{eofer-14}
\liminf_{n \to \infty} \PPP \times
\PPP' \left( |I_{n, \de}| \geq n^{1-a-\ep} \right) \geq 1/2
\end{equation}
and
\beq{eofer-99}
\lim_{n \to \infty}
\PPP \times \PPP' ( \inf_{x \in I_{n, \de} } L^{X,X'}_{n} (x)
\leq \frac{(1-\ep)^2 a^2}{4 \pi^2} (\log n)^4 ) = 0\,.
\end{equation}
\end{lemma}

\noindent
{\bf Proof of Lemma \ref{eSRW-occ}:} We first prove (\ref{eofer-99}).
Fixing $n$, let
$\si_z$ denote the time it takes $\{ X_i \}$ to
complete $(1-2\de) m_{k(n)}$ excursions between $\partial D_z(r_n)$ and
$\partial D_z( R_n)$, after first hitting $\partial D_z(R_n)$, and
denote by $Y_j(z)$ the occupation measure of $z$ by $\{ X_i \}$
during its $j$-th such excursion. Thus,
fixing the starting point of the $j$-th excursion,
$Y_j(z)$ is distributed like
$L_{T_{R_n}}(0)$
where $T_{R_n}$ is the first hitting time of $\partial D_0(R_n)$
(when starting at the corresponding $z' \in \partial D_0(r_n)$).
We note that
\begin{equation}
\E^0\(\lc L_{T_{R_n}}(
0)\rc^k\)=k!\lc
G_{R_n}(0,0)\rc^k\hspace{.3in}k=1,2,\ldots\label{ej.30}
\end{equation}
where $G_{R_n}(x,y)$ is Green's function for $D_0(R_n)$. This is a
simple case of Kac's
moment formula,
see
\cite{FP}. We will only need $k=1$ and $2$. We will use the
abbreviation $G_n=G_{R_n}(0,0)$.

Set $\la=\de^{1/4}$, $\ell_n=(1-2\la)^2/( \de \log n )$.  As
$$
L_{\si_z}(z) \geq \sum_{j=1}^{(1-2\de) m_{k(n)}} Y_j(z) \,,
$$
it follows by the strong Markov property of $\{ X_i \}$ at the start
of each of these excursions that
\beqn{eofer-8}
&&
\sup_{y \in \Z^2} \PPP^y \( {L_{\si_z}(z)\over G_n } \leq (1-2\de) m_{k(n)}
\ell_n \)\\
 && \leq
e^{\la (1-2\de) m_{k(n)} \ell_n} \sup_{y \in \Z^2} \E^y ( e^{-\la
{L_{\si_z}(z)\over G_n } } )
\nn
\\
&& \leq
\Big[ e^{\la \ell_n} \sup_{z' \in \partial D_z(r_n)}
 \E^{z'} ( e^{- \la L_{T_{R_n}}(z)/ G_n } ) \Big]^{(1-2\de) m_{k(n)}}
\nn
\\
&& = \Big[ e^{\la \ell_n} \sup_{z \in \partial D_0(r_n)}
 \E^z ( e^{-\la L_{T_{R_n}}(0)/ G_n } ) \Big]^{(1-2\de) m_{k(n)}}
\;,\nn
\eeqn
with the last identity following from
the translation invariance of
SRW. We first study the
quantity appearing in the the first line of (\ref{eofer-8}). By \cite[Theorem
1.6.6]{lawler}, for all $n$ large enough,
\beq{eofer-9}
G_{R_n}(0,0)
= \frac{2}{\pi} \log R_n+O(1) =\de \frac{2}{\pi} \log n+O(1).
\end{equation}
Hence for $n$ large enough,
\begin{eqnarray}
(1-2\de) m_{k(n)} \ell_n G_n &\geq&
(1-2\de) (1-2\de^{1/4})^3 \frac{2}{ \pi} a k^2(n)\label{eofer-13}\\
&\geq& (1-2\de^{1/4})^5 \frac{a}{2 \pi} (\log n)^2.\nn
\end{eqnarray}

If $z \in D_0(\sqrt{n})$ is $n, \de$-admissible,
then necessarily
\beq{eofer-15}
L^{X,X'}_{n}(z) \geq  L_{\si_z}(z)  L'_{\si'_z}(z)
\;,
\end{equation}
where $\si'_z$ denotes the time it takes $\{ X'_i\}$ to
complete $(1-2\de) m_k(n)$ excursions between $\partial D_z(r_n)$ and
$\partial D_z( R_n)$, after first hitting $\partial D_z(R_n)$.

There are at most $\pi n$ lattice points
$z \in D_0(\sqrt{n})$, hence by  \req{eofer-15} and \req{eofer-13}
we conclude that for $\de \leq \de_0$ and all
$n$ large enough
\begin{eqnarray}
&& \PPP \times \PPP' ( \inf_{z \in I_{n,\de} } L^{X,X'}_{n}(z) \leq
 \frac{ (1-2\de^{1/4})^{10} a^2}{4 \pi^2} (\log n)^4 ) \label{ej.32}\\
&&   \leq
2 \pi n
\sup_{y \in \Z^2} \PPP^y (  L_{\si_z}(z) \leq (1-2\de)
m_{k(n)} \ell_n G_n)\nn
\end{eqnarray}
and \req{eofer-99} will follow once we show that the right hand side of
\req{eofer-8}
is bounded by $n^{-2}$, which we now do.

Let $T_0$ be the first
hitting time of $0$ and set $q_{n,z}:=\PPP^z ( T_0<T_{R_n} )$.
 The Markov property followed by the inequality $e^{-t}\leq
1-t+t^2/2\,;\,t\geq 0$
shows that
\begin{eqnarray} \hspace{.3in}
&&
\E^z ( e^{-\la L_{T_{R_n}}(0)
/G_n
} )\label{emoment1}\\
&& =\PPP^z ( T_0>T_{R_n} )+
\PPP^z ( T_0<T_{R_n} )\E^0( e^{-\la  L_{T_{R_n}}(0)/ G_n})
\nn\\
& &\leq  1-q_{n,z} +q_{n,z}\(1-\la\E^0\({L_{T_{R_n}}(0)\over G_n}\)
+{\la^2\over 2}\E^0\(\lc {L_{T_{R_n}}(0)\over G_n}\rc^2\)\)
\nn
\end{eqnarray}
and by (\ref{ej.30}) this gives
\begin{equation}
\E^z ( e^{-\la L_{T_{R_n}}(0)} )\leq 1-q_{n,z}(\la-\la^2)
\leq e^{- q_{n,z} (\la - \la^2)}.\label{emoment1a}
\end{equation}

By \cite[Proposition 1.6.7]{lawler}, for all $n$ large enough, and all $z \in
\partial D_0(r_n)$
\begin{equation}
q_{n,z}
= \frac{\log (R_n/r_n)}{\log R_n}+O(\frac{1}{\log^2 R_n})
\geq (1-2\la)/( \de \log n )\label{eofer-10}
\end{equation}
Then with $q_n := (1-2\la)/( \de \log n )$
 we see from (\ref{emoment1a}) and (\ref{eofer-10}) that
\begin{equation}
\sup_{z \in \partial D_0(r_n)} \E^z ( e^{-\la L_{T_{R_n}}(0)
/G_n
} )
\leq e^{- q_n (\la - \la^2)}.
\end{equation}
Since $\ell_n=(1-2\la)q_n$
we deduce that
\beq{eofer-11}
e^{\la \ell_n}
\sup_{z \in \partial D_0(r_n)} \E^z ( e^{-\la L_{T_{R_n}}(0)
/G_n
} ) \leq
e^{-\la^2 q_n}
\end{equation}
and for all $\de \leq \de_0(a)$ and $n$ large enough
$$
e^{-\la^2 q_n (1-2\de) m_{k(n)}}
\leq e^{-a k^2(n)/(\la \log n)} \leq n^{-2} \;.
$$
The last two displays show that the right hand side of \req{eofer-8}
is bounded by $n^{-2}$.
As mentioned after (\ref{ej.32}) this completes the proof of \req{eofer-99}.

We next turn to the proof of (\ref{eofer-14}). It is here that we use strong
approximation, so we need to introduce further notation concerning Brownian
paths.
Let $M_k^x$ denote the number of
excursions of $\{W_t\}$ from $\partial D(x,e^{-k+1})$ to
$\partial D(x,e^{-k})$ prior to $\bart$ and
${M_k'}^x$ denote the corresponding number of excursions for
the process $\{W'_t\}$.  Fix $\de,\,\bb>0$ and let $\UU_k(\bb)$ be a fixed
maximal
collection of points $x_j\in D(0,1)\cap D^c(0,\bb)$ with $|x_i-x_j|\geq
\de^2e^{-k}$ for all
$x_i,x_j\in
\UU_k(\bb)$. We say that a point $x \in \UU_k(\bb)$ is {\bf  $k$-admissible}
if $M_k^x \wedge {M'_k}^x \geq (1-2\de) m_k$ and denote by
$\wh{\UU}_k(\bb)$ the set of $k$-admissible points.

We now show how to derive (\ref{eofer-14}) from the following lemma whose
proof is
momentarily deferred.
\begin{lemma}
\label{ekey-BM}
For $\bb>0$ sufficiently small
\begin{equation}
\liminf_{k \to \infty}\PPP \times \PPP' \(|\wh{\UU}_k(\bb)|\geq
e^{2(1-a-2\de)k}\)\geq
7/8\label{ej.2}
\end{equation}
\end{lemma}

Using (\ref{ej.2}) together with the fact that $\bart \vee \bart'<\ff$ a.s.
we can
find $K<\infty$ so that
\begin{equation}
\liminf_{k \to \infty}
\PPP \times \PPP' \(|\wh{\UU}_k(\bb)|\geq e^{2(1-a-2\de)k}
;\,\, \bart
\vee \bart'
\leq K\)\geq 3/4.\label{ej.4}
\end{equation}

By Brownian scaling and the
multidimensional strong approximation of \cite[Theorem 1]{Ei} we
may
construct for each $n$ independent SRW's $\{X_i\}$, $\{X'_i\}$
and independent Brownian motions $\{W_t : t \in [0,K] \}$,
$\{ W'_t : t \in [0,K]\}$ on the same probability space such that
$\PPP \times \PPP' ( \BB_n ) \to 1$ as $n \to \infty$, where $\BB_n$ is the set
\[\{ \sup_{0 \leq t \leq K}
|W_t - \frac{\sqrt{2K}}{\sqrt{n}} X_{[t n/K]}| \leq \de e^{-k(n)},
\sup_{0 \leq t \leq K}
|W'_t - \frac{\sqrt{2K}}{\sqrt{n}} X'_{[tn/K]}| \leq \de e^{-k(n)} \}
\,.\]
Recall that if $\bart \vee \bart' \leq K$ then the event that
$x \in D(0,1)$ is $k$-admissible is measurable on
$\sigma(W_t, W'_t : t \leq K)$. Hence, for all $n$ large enough,
if in addition $\BB_n$ holds,
then to each $k(n)$-admissible $x \in D(0,1)$
corresponds $\widetilde x \in \Z^2$ nearest to $\sqrt{n/(2K)} x$
that is $n, \de$-admissible.
Since
$\PPP \times \PPP' ( \BB_n ) \to 1$ as $n \to \infty$, (\ref{eofer-14})
follows.\qed

\noindent
{\bf Proof of Lemma \ref{ekey-BM}:} We begin by using the techniques of
\cite{DPRZ4}
to find many points in whose neighborhood both Brownian paths have large
occupation measure.  Let
\beq{eCa-def}
\CC_a:=\Big\{x\in D(0,1): \lim_{\ep\to 0}
\frac{\mu_{\bart}^{w} (D(x,\eps))}{\eps^2 \left(\log \ep\right)^2}= a\Big\},
\end{equation}
and let
$\CC'_a$ denote the corresponding set for
the process $\{W'_t\}$.
\begin{lemma}\label{lem-codim}
\begin{equation}
\PPP \times \PPP' \lc \dim (\CC_a\cap \CC'_a)=2-2a\rc=1.\label{ej.1}
\end{equation}
\end{lemma}

\noindent
{\bf Proof of Lemma \ref{lem-codim}}. We follow the
proof of \cite[Theorem 1.3]{DPRZ4},
but we now say that the indicator function $Y(n,i)$
introduced in \cite[Section 3]{DPRZ4} is equal to $1$
iff $x_{n,i}$ is n-perfect both for $W$ and
$W'$. By
independence, the bounds for the first moment and covariance of $Y(n,i)$
which appear
in \cite[Lemma 3.2]{DPRZ4} now have $a$ replaced by $2a$. The rest of the
proof now
proceeds exactly as in \cite{DPRZ4}.
\qed

The following lemma allows us to obtain large excursion
counts from the large occupation times provided by the previous lemma.
We abbreviate $\rho_k=e^{-k}$.
\begin{lemma}\label{lem-ttoe}
For any
$\de,\bb,
\ga>0$ we can find $k_0=k_0(\om)<\ff$ a.s. such that for all
$k>k_0$ and $x\in \UU_k(\bb)$,  if
\begin{equation}
\mu^{W}_{\bart} (D(x,\rho_k)) \geq
a(1-\ga^2) |\rho_k\log \rho_k|^2 \label{etoe.1}
\end{equation}
then $M_k^x\geq (1-\ga)m_k$.
\end{lemma}

\noindent
{\bf Proof of Lemma \ref{lem-ttoe}}.
We will say that $x\in \UU_k(\bb)$ is *-thick if it satisfies (\ref{etoe.1}).
Assuming $k$ is large enough so that $0 \notin D(x,\rho_k)$ for all
$x \in \UU_k(\bb)$,
we let $\tau_{l,k}$
denote the occupation measure of $D(x,\rho_k)$ during the
$l$-th excursion of $W$ between $\partial D(x,\rho_k)$ and
$\partial D(x,\rho_{k-1})$.
Then, with $\E(\tau_{l,k})=\rho_k^2$ and
$m_k'=(1-\ga) m_k = a(1-\ga) (\log \rho_k)^2$,
we have that for some universal constant $C<\infty$ and all
$x \in \UU_k(\bb)$,
\beaa
\PPP ( M_k^x \leq m_k' \,,\,\mbox{$x$ is *-thick}) &\leq&
\PPP ( \sum_{l=1}^{m_k'} \tau_{l,k} \geq
a(1-\ga^2) |\rho_k\log \rho_k|^2 ) \\
&=&
\PPP ( \frac{1}{m_k'} \sum_{l=1}^{m_k'} \tilde{\tau}_{l,k} \geq  \ga)
\leq e^{-\ga^2 m_k/C} \;,
\eeaa
where $\tilde{\tau}_{l,k}:= \tau_{l,k}/\E (\tau_{l,k}) -1$ ,
and the last inequality follows by the methods used in the proof of
\cite[Lemma 6.4]{DPRZ4}. Consequently,
\beq{eofer-16}
 \sum_{k=1}^\ff\sum_{x \in \UU_k(\bb)}\PPP ( M_k^x \leq  m_k'\,,\,\mbox{$x$ is
*-thick} ) \leq
 \sum_{k=1}^\ff |\UU_k(\bb)| e^{-\ga^2 m_k/C}<\ff.
\end{equation}
\qed

We can now complete the proof of  Lemma \ref{ekey-BM}. It follows from
(\ref{ej.1})
that for some $\bb>0$ sufficiently small and all
$k$ sufficiently large, with probability $\geq 7/8$ the
set of
$x\in D(0,1)\cap
D^c(0,2\bb)$ with
\begin{equation}
\min
\(\frac{\mu_{\bart}^{w} (D(x,(1-\de^2)\rho_k))}{\rho_k^2 \left(\log
\rho_k\right)^2}\,,\frac{\mu_{\bart'}^{w'} (D(x,(1-\de^2)\rho_k))}{\rho_k^2
\left(\log
\rho_k\right)^2}\,\)\geq (1-2\de^2)a\label{ej.3}
\end{equation}
has Hausdorff dimension $\geq 2-2a-\de$. Let $\wt{\UU}_k$ be the set of
points in
$\UU_k(\bb)$
which are within $\de^2\rho_k$ of the set in (\ref{ej.3}).
Using Lemma
\ref{lem-ttoe} it is easy to check that each point in $\wt{\UU}_k$ is
$k$-admissible.
 Since $\{ D(x,\de\rho_k) : x \in \wt{\UU}_k \}$ is a cover
of the set in (\ref{ej.3}) by sets of maximal diameter $\de^2\rho_k$, it
follows
that
\beq{eofer-17}
\liminf_{k \to \infty} |\wt{\UU}_k| (\de^2 \rho_k)^{2-2a-2\de}= \infty .
\end{equation}
Our lemma now follows.
\qed

\section{Complements and unsolved problems}\label{sec-comp}

\noindent
$\Large \bullet \;$
By Brownian scaling, for any deterministic $0<r<\infty$,
the set $D(0,1)$ and $\bart,\bart'$ can be replaced by
$D(0,r)$  and $\bart_r = \inf\{s:\,|W_s| = r\},\bart'_r = \inf\{t:\,|W'_t|
= r\}$, without changing the conclusion
of Theorem \ref{theo-1}. Similarly, one may replace $\II_{\bart,\,\bart'}$ by
$\II_{S,\,T}$ in this theorem, for any deterministic $0<S,T<\infty$.
Moreover, from its proof we have that 
\req{01} remains valid when the limit in $\eps$ is replaced by
$\liminf$ or $\limsup$ and when considering the 
set of points $x$ for which this limit ($\liminf$, $\limsup$, 
respectively), is at least $a^2$. 

\medskip
\noindent $\Large \bullet\;$
Next, we discuss briefly the packing dimension analogue of Theorem
\ref{theo-1}; consult Mattila (1995) for background on packing dimension,
Minkowski dimension and their relation.
The set of consistently thick points ${\mbox{\sf CThickInt}}_{\ge a^2}$,
 defined in (\ref{defcthi}), has different
packing dimension from the set ${\mbox{\sf ThickInt}}_{\ge a^2}$, defined in
(\ref{defthi}). Namely, for every $a \in (0,1]$,
\begin{equation} \label{packcthi}
\dim_P({\mbox{\sf CThickInt}}_{\ge a^2}) =2-2a \, , \;\; a.s.
\end{equation}
\begin{equation} \label{packthi}
\dim_P({\mbox{\sf ThickInt}}_{\ge a^2}) =2\, \;\; a.s.
\end{equation}
To justify (\ref{packcthi}), we use the notation of Section 2.
The sets $\AA_n$, defined in (\ref{defaa}),
satisfy
\begin{equation} \label{eqdelta}
|\AA_n|\leq (\tilde{\ep}_n)^{(1-11\de)2 a  -2}
\end{equation}
for all large $n$, by (\ref{3.3jj}) and Borel-Cantelli.

Recall the discs $\VV_{n,j}=D(x_j,\de\tilde{\ep}_n)$
defined after (\ref{3.3jj}), and denote
$\VV_n=\cup_{j \in \AA_n} \VV_{n,j}$. By (\ref{eqdelta}),
the upper Minkowski
dimension of $\VV_{\ell}^*= \cap_{n \ge \ell}  \VV_n$ is at most
$2-(1-11\de) 2a $.
It is easy to see that ${\mbox{\sf CThickInt}}_{\ge a^2} \subset \cup_{\ell
\ge 1}
\VV_{\ell}^*$, whence $\dim_P({\mbox{\sf CThickInt}}_{\ge a^2}) \le
2-(1-11\de) 2a $.
Since $\delta$ can be taken arbitrarily small,
while $\dim_P({\mbox{\sf CThickInt}}_{\ge a^2})\geq \dim({\mbox{\sf
CThickInt}}_{\ge
a^2})$, this proves (\ref{packcthi}).

\smallskip

To prove (\ref{packthi}), it clearly suffices
to consider $a=1$. Recall that
  $\bart = \inf\{t:\,|W_t| = 1\}$. For each $n \ge 1$,
 let
$$
V_n:= \bigcup_{0<\eps< 1/n} \Big \{0 < t< \bart:\;
\frac{\II_{\bart,\,\bart'} (D(W_t,\eps))}{\eps^2
\left(\log \frac{1}{\eps}\right)^4} > 1-1/n \Big\}.
$$
The sets $V_n$ are open by the path continuity of Brownian motion.
Moreover, it
is easy to check, e.g. by applying Theorem \ref{theo-2} with an arbitrary
$T$ replacing $\bart$ there, and using the shift invariance of Brownian motion,
that for any $n \ge 1$,
almost surely
$V_n$ is a dense subset of $(0,\bart)$;
by  \cite[Corollary 
2.4, part (i)]{DPRZ1}, $\dim_P(\cap_n V_n)=1$ a.s.
The set
\begin{equation}\label{thicktime}
\Big\{0\leq t\leq \bart:\;
\limsup_{\eps\to 0} \frac{\II_{\bart,\,\bart'} (D(W_t,\eps))}{\eps^2
\left(\log \frac{1}{\eps}\right)^4} \ge 1\Big\},
\end{equation}
contains $\cap_n V_n$, so it has packing dimension 1.
Finally, ${\mbox{\sf ThickInt}}_{\ge 1}$
is the image under planar Brownian motion of the set in (\ref{thicktime});
hence the uniform doubling of packing dimension by
 planar Brownian motion, see \cite[Corollary 5.8]{Perkins-Taylor},
yields (\ref{packthi}).

\medskip
\noindent $\Large \bullet\;$
The situation for intersections of $m$ independent planar Brownian motions is
completely analogous to that of Theorems \ref{theo-2} and \ref{theo-1}.
Specifically, let $W^{(1)}_{s_1},\ldots, W^{(m)}_{s_m}$ denote $m $
independent planar Brownian motions, and define the $m$-fold projected
intersection
local time by
\beaa
&&
\II_{m,S_1,\ldots, S_m}(A)
\\
&&=\lim_{\ep\rar
0}\,\,\pi^{m-1}\int_0^{S_1}\cdots\int_0^{S_m}1_{A}(W^{(1)}_{s_1})
\prod_{j=2}^m f_\ep(W^{(1)}_{s_1}-W^{(j)}_{s_j})\,ds_1\ldots\,ds_m\nn
\eeaa
where $f_\ep$ is any approximate $\de$-function. It is known that the limit
(\ref{i1.1}) exists a.s. and in all $L^p$ spaces, and that
$\II_{m,S_1,\ldots, S_m}(\cdot)$ is a
measure  supported on \[\{x\in\reals^2|
\,x=W^{(1)}_{s_1}=\cdots=W^{(m)}_{s_m}
\,\,\mbox{\rm  for some
$0\leq s_1\leq S_1,\ldots,\,0\leq s_m\leq S_m$}\},\] see
\cite[Theorem 1, Chapter VIII]{Le Gall}.
 Note that
$\II_{1,\bart}$ is the occupation measure studied in \cite{DPRZ4}. Let
$\bart^{(j)} =
\inf\{s:\,|W^{(j)} _s| = 1\}$. Then,
for any $m\geq 1$
\begin{equation}
\lim_{\eps\to 0}
\sup_{x \in \reals^2} \frac{\II_{m,\bart^{(1)},\ldots, \bart^{(m)}}
(D(x,\eps))}{\eps^2
\left(\log \frac{1}{\eps}\right)^{2m}}=(2/m)^m\,,
\hspace{.6in}a.s.
\end{equation}
and for any $0<a\leq 2/m$,
\beq{01m}
\dim\Big\{x\in D(0,1):\;
\lim_{\eps\to0} \frac{\II_{m,\bart^{(1)},\ldots, \bart^{(m)}}
(D(x,\eps))}{\eps^2
\left(\log \frac{1}{\eps}\right)^{2m}}= a^m\Big\}= 2-ma
\hspace{.3in}a.s.
\end{equation}

\medskip
\noindent $\Large \bullet\;$
The Hausdorff dimension of the set of thick intersection points
for two independent Brownian motions in $\reals^3$
was recently determined by K\"{o}nig and M\"{o}rters.
Specifically, let $B(x,r)$ denote the ball of radius $r$ centered 
at $x \in \reals^3$ and  $\JJ(B(x,r))$ the total intersection local
time in $B(x,r)$ for two independent Brownian motions in $\reals^3$ 
(see \cite[Section 2.1]{Mor} for various equivalent definitions
of $\JJ(\cdot)$). In \cite[Theorem 1.4]{Mor}, 
K\"{o}nig and M\"{o}rters show that for any $0<a\leq \rho^*/2$, 
almost surely,
\beq{conj3}
\dim\Big\{x\in R^3:\;
\limsup_{\eps\to0} \frac{\JJ (B(x,\eps))}{\eps
\left(\log \frac{1}{\eps}\right)^2}= a^2\Big\}= 1-2a/\rho^*\,,
\end{equation}
where the non-random $\rho^*>0$ is the solution of the explicit
variational formula \cite[(1.9)]{Mor}. The analog of \req{conj3} 
for consistently thick points, that is with  
$\liminf$ instead of $\limsup$, involves a different gauge 
function and remains an open problem.

\medskip
\noindent $\Large \bullet\;$
In \cite{DPRZ2} we analyzed `thin points' for the Brownian occupation
measure, establishing that
\begin{equation}
\lim_{\eps\to 0}\inf_{t\in [0,1]} \frac{\mu_\bart^W (D(W_t,\eps))}{\eps^2
/\log \frac{1}{\eps}}=1\,,
\hspace{.6in}a.s.\label{thin.1}
\end{equation}
with the multi-fractal spectrum
\begin{equation}
\dim\{x\in D(0,1):\;
\liminf_{\eps\to 0} \frac{\mu_\bart^W (D(x,\eps))}{\eps^2
/\log \frac{1}{\eps}}= a\}= 2-2/a
\hspace{.3in}a.s.\label{thin.2}
\end{equation}
for any fixed $a>1$. In the present paper, we analyzed
Brownian projected intersection local time
where it is exceptionally `thick'.
The analysis of the corresponding `thin intersection points'
is needed to describe completely the multi-fractal structure
of this measure, and remains an open problem.

\medskip
\noindent
{\bf Acknowledgment}
We thank Wolfgang K\"{o}nig and Peter M\"{o}rters for
providing us with the preprint \cite{Mor}.

\bigskip
\noindent
\begin{tabular}{lll}  & Amir Dembo& Yuval Peres\\
  & Departments of
Mathematics &  Dept.\ of Statistics, UC Berkeley\\
& and of Statistics & Berkeley, CA 94720  and \\
 &Stanford University & Institute of Mathematics\\
 &Stanford, CA 94305 & Hebrew University, Jerusalem, Israel \\
&amir@math.stanford.edu&peres@stat.berkeley.edu\\
& &\\
& & \\
& & \\
  & Jay Rosen & Ofer Zeitouni\\
  & Department of
Mathematics& Department of Electrical Engineering\\
 &College of Staten Island, CUNY& Technion\\
 &Staten Island, NY 10314& Haifa 32000, Israel\\
&jrosen3@idt.net& zeitouni@ee.technion.ac.il
\end{tabular}

\end{document}